\documentclass{amsart} %12pt
\usepackage{enumerate}
\usepackage{amsmath,amssymb,amsxtra,amsbsy,amsthm}
\usepackage[mathscr]{eucal}
\usepackage[all]{xy}

\newtheorem{thm}{Theorem}[section]
\newtheorem{lemma}[thm]{Lemma}
\newtheorem{prop}[thm]{Proposition}
\newtheorem{cor}[thm]{Corollary}
\theoremstyle{definition}

\theoremstyle{remark}
\newtheorem{rem}[thm]{Remark}
\newtheorem{rems}[thm]{Remarks}

\numberwithin{equation}{section}

\newcommand{\vanish}[1]{\relax}       % Auskommentieren von Text
\newcommand{\Ell}[2][\relax]{%
   \ifx#1\relax \mathrm{L}^{\mathrm{#2}}
   \else \mathrm{L}^{\mathrm{#2}}_{\mathrm{#1}}
   \fi}
\newcommand{\loc}{\text{\upshape{\tiny loc}}}

\DeclareMathOperator{\rmT}{T}
\newcommand{\ud}{\mathrm{d}}
\newcommand{\pfeil}{\rightarrow}
\newcommand{\tpfeil}{\mapsto}
\DeclareMathOperator{\dom}{dom}
\DeclareMathOperator{\ran}{ran}
\DeclareMathOperator{\Lin}{\mathcal{L}}
\newcommand{\T}{\mathbb{T}}
\newcommand{\R}{\mathbb{R}}
\newcommand{\N}{\mathbb{N}}
\newcommand{\Z}{\mathbb{Z}}

\newcommand{\C}{\mathbb{C}}
\newcommand{\D}{\mathbb{D}}

\newcommand{\Sum}[2][\relax]{%
 \ifx#1\relax \sideset{}{_{#2}}\sum
 \else \sideset{}{^{#1}_{#2}}\sum
 \fi}                                   % Summen mit Index nebengestellt
\newcommand{\Lap}{\mathcal{L}}
\newcommand{\Fourier}{\calF}
\newcommand{\fourier}[1]{\widehat{#1}}

\newcommand{\norm}[2][\relax]{%                             % norm
   \ifx#1\relax \ensuremath{\left\Vert#2\right\Vert}
   \else \ensuremath{\left\Vert#2\right\Vert_{#1}}
   \fi}
\newcommand{\calA}{\mathcal{A}}
\newcommand{\abs}[1]{\left\vert#1\right\vert}   % Absolutbetrag math. Modus

\newcommand{\Ha}{\mathrm{H}}
\DeclareMathOperator{\re}{Re}          % Realteil
\DeclareMathOperator{\im}{Im}          % Imagin\"arteil
\DeclareMathOperator{\sgn}{sgn}        % Signum
\newcommand{\suchthat}{\,\,|\,\,}     % Durch Aussonderung definierte Mengen
\newcommand{\st}{\,\,|\,\,}     % Durch Aussonderung
\newcommand{\Be}{\mathrm{B}}
\newcommand{\eM}{\mathrm{M}}
\newcommand{\BV}{\mathrm{BV}}
\newcommand{\car}{\mathbf{1}}           % Konstante Funktion 1
\DeclareMathOperator{\supp}{supp}       % topologischer Tr\"ager
\newcommand{\calF}{\mathcal{F}}

\newcommand{\vphi}{\varphi}           % Abk\"urzung
\newcommand{\vpsi}{\psi}           % Abk\"urzung
\newcommand{\ohne}{\setminus}

\newcommand{\sprod}[2]{\ensuremath{\left<#1,#2\right>}}
\newcommand{\nach}{\circ}
\newcommand{\ue}{\mathrm{e}}

\newcommand{\res}[1]{ \big|_{#1}}      % Einschr\"ankung
\newcommand{\konj}[1]{\overline{#1}}   % Konjugation
\newcommand{\mnorm}[1]{\|#1\|}         % norm im Fliesstext
\DeclareMathOperator{\card}{card}  % Anzahl, alternativ: \#
\newcommand{\calM}{\mathcal{M}}
\newcommand{\Ce}{\mathrm{C}}

\newcommand{\Exp}{\mathbb{E}}      % Expectation
\newcommand{\tensor}{\otimes}
\newcommand{\ce}{\mathrm{c}}

\newcommand{\Wee}[2][\relax]{%
   \ifx#1\relax \mathrm{W}^{\mathrm{#2}}
   \else \mathrm{W}^{\mathrm{#2}}_{\mathrm{#1}}
   \fi}
\newcommand{\calT}{\mathcal{T}}
\newcommand{\floor}[1]{\lfloor#1\rfloor}

\newcommand{\AM}{\mathcal{AM}}
\newcommand{\sector}[1]{\Sigma_{#1}}
\newcommand{\PV}{\mathrm{PV}}

\newcommand{\Mlt}{\mathcal{M}}
\newcommand{\rmU}{\mathrm{U}}
\newcommand{\rmA}{\mathrm{A}}
\newcommand{\UMD}{{\scshape umd}}

\newcounter{aufzi}

\newcounter{aufzii}

\newcounter{aufziii}
\newenvironment{aufziii}{\begin{list}{ {\upshape\arabic{aufziii})}}{
        \usecounter{aufziii}
        \topsep1ex
        \parsep0cm
        \itemsep1ex
        \leftmargin0.8cm
        \labelwidth0.5cm
        \labelsep0.3cm
}}
{\end{list}}

\newcounter{aufziv}

\begin{document}

\title[Transference Principles for Semigroups]{Transference Principles for Semigroups and a Theorem of Peller}

\author{Markus Haase}%
\address{Delft Institute of Applied Mathematics\\ Delft University
of Technology\\
P.O.Box 5031\\ 2600 GA Delft\\ The Netherlands%\\[1em]
}%

\email{m.h.a.haase@tudelft.nl}%

\begin{abstract}
A general approach to transference principles for discrete and continuous
operator (semi)groups is described. This  allows to recover
the classical transference results of 
Calder\'on, Coifman and Weiss and of Berkson, Gillespie and Muhly 
and  the more recent one of the author. 
The method is applied to derive a new transference
principle for (discrete and continuous) operator
semigroups that need not be groups. 
As an application,
functional calculus estimates
for bounded operators with 
at most polynomially growing powers are derived, culminating in 
a new proof of classical results by Peller from 1982.  
The method allows a generalization of his results away from
Hilbert spaces to $\Ell{p}$-spaces and
--- involving the concept of $\gamma$-boundedness --- to general Banach
spaces. Analogous results for strongly-continuous one-parameter
(semi)groups are presented as well. Finally, an application is given to
singular integrals for one-parameter semigroups.
\end{abstract}

\date{\today}



% \linenumbers


\renewcommand{\subjclassname}{\textup{2000} Mathematics Subject Classification}

\maketitle

\section{Introduction and Summary}\label{s.intro}

The purpose of this article is twofold. The shorter part 
(Section \ref{s.ti}) is
devoted to a  generalization of the classical 
transference principle of Calder\'on, Coifman and Weiss. 
In the  major part
(Sections \ref{s.grp}--\ref{s.si}) we give applications of this new abstract results
to discrete and continuous operator (semi)groups; in particular
we shall recover and generalize important results of Peller \cite{Pel82}.

\medskip

In the classical transference principle(s) 
the objects under investigation 
are derived operators of the form
\begin{equation}\label{intro.e.grp}
  \rmT_\mu := \int_G T(s)\, \mu(\ud{s})
\end{equation}
where $G$ is a locally compact group and $\rmT = (T(s))_{s\in G} : G \pfeil \Lin(X)$
is a bounded strongly continuous representation of $G$ on a Banach space $X$. 
The integral \eqref{intro.e.grp} has to be understood in the strong sense, i.e.,
\[ \rmT_\mu x = \int_G T(s)\, \mu(\ud{s}) \qquad (x\in X).
\]
Since such operators occur in a variety of situations, the applications
of transference principles are manifold, and
the literature on this topic is vast. We therefore restrict 
ourselves to mentioning only a few 'landmarks' which  we regard as
most important for the understanding of the present paper.

Originally, Calder\'on \cite{Cal68}  considered representations on $\Ell{p}$
induced by  a $G$-flow of measure-preserving transformations of the underlying
measure space. His considerations were motivated by ergodic theory and
his aim was to obtain maximal inequalities. 
Subsequently, Coifman and Weiss \cite{CoiWei,CoiWei77}
shifted the focus to norm estimates and were able to drop Calder\'on's assumption
of an underlying measure-preserving $G$-flow towards general $G$-representations
on $\Ell{p}$-spaces. Some years later, 
Berkson, Gillespie and Muhly \cite{BerGilMuh89b} were able to generalize the method
towards general Banach spaces $X$. 
However, the representations considered in these works were still (uniformly) bounded.
In the continuous one-parameter case (i.e., $G = \R$)
Blower \cite{Blo00} showed 
that the original proof method could fruitfully
be applied also to non-bounded representations.
However, his result was in a sense 'local' and did not take
into account the growth rate of the group $(T(s))_{s\in \R}$ at infinity.
In \cite{Haa07b} we re-discovered Blower's result and in \cite{Haa09b}
we could refine it towards a 'global' transference result for strongly continuous
one-parameter groups, cf.~also Section \ref{s.grp} below.

\medskip

In the present paper, more precisely in Section \ref{s.ti}, we develop
a {\em method of generating transference results} and show in Section \ref{s.grp}
that the known transference principles
(the classical Berkson-Gillespie-Muhly result and the central results
of \cite{Haa09b}) are special instances of it.
Our method has three important new features. Firstly, it allows to
pass from {\em groups} (until now the standard assumption) to {\em semigroups}.
More precisely, we consider closed sub-semigroups $S$ of a locally compact group $G$
together with a strongly continuous representation $\rmT: S \pfeil \Lin(X)$
on a Banach space, and try to estimate the norms of operators of the form
\begin{equation}\label{intro.e.sgp}
\rmT_\mu =  \int_S T(s) \, \mu(\ud{s})
\end{equation}
by means of the transference method. The second feature is the role
of {\em weights} in the transference procedure, somehow hidden in the classical
version.  Thirdly, our account brings to light
the {\em formal structure} of the transference argument: in a first step
one establishes a {\em factorization} of the operator \eqref{intro.e.sgp}
over a convolution  (i.e., Fourier multiplier) operator 
on a space of $X$-valued functions on $G$; then,
in a second step, one uses this factorization to estimate the norms;
and finally, one may vary the parameters to optimize the obtained
inequalities. So one can briefly subsume our method under the scheme
\begin{quote}
\begin{center}
factorize --- estimate --- optimize,
\end{center}
\end{quote}
where we use one particular way of constructing the initial factorization.
One reason for the power of the method lies 
in choosing different weights in the factorization, allowing for the
optimization  in the last step. The second reason lies in
the purely formal nature of the factorization; this allows to
 re-interpret the same factorization involving different function spaces.

\medskip

The second part of the paper (Sections \ref{s.ops}--\ref{s.si})
is devoted to applications of the 
transference method. These applications deal exclusively with the 
cases $S = \Z, \Z_+$ and 
$S= \R, \R_+$,
which we for 
short call the {\em discrete} and the {\em continuous} case, respectively. 
However, let us 
point out that the general transference method of Section \ref{s.ti}
works even for sub-semigroups of {\em non-abelian} groups.  
 
To clarify what kind of applications we have in mind, let us
look at the discrete case first. Here the semigroup consists of the powers
$(T^n)_{n \in \N_0}$ of one single bounded operator $T$, and 
the derived operators \eqref{intro.e.sgp} take the form
\[ \Sum{n \ge 0} \alpha_n T^n
\]
for some (complex) scalar sequence $\alpha = (\alpha_n)_{n \ge 0}$. In order to avoid
convergence questions, we suppose that $\alpha$ is a finite sequence,
hence 
\[ \fourier{\alpha}(z) :=  \Sum{n \ge 0} \alpha_n z^n
\]
is a complex polynomial. One usually writes
\[  \fourier{\alpha}(T) := \Sum{n \ge 0} \alpha_n T^n 
\]
and is interested in continuity properties of the {\em functional calculus}
\[ \C[z]  \pfeil \Lin(X),\qquad f \tpfeil f(T).
\]
That is, one looks for a function algebra norm $\norm{\cdot}_\calA$ on $\C[z]$
that allows
an estimate of the form
\begin{equation}\label{intro.e.fc-est}
 \norm{f(T)} \lesssim \norm{f}_\calA \qquad (f\in \C[z]).
\end{equation}
A rather trivial instance of \eqref{intro.e.fc-est}  is based on the estimate
\[ \norm{f(T)} = \norm{\Sum{n\ge 0} \alpha_n T^n} 
\le \Sum{n\ge 0}  \abs{\alpha_n} \norm{T^n}.
\]
Definining the positive sequence $\omega= (\omega_n)_n$ 
by $\omega_n := \norm{T^n}$, we hence have 
\begin{equation}\label{intro.e.fc-triv}
 \norm{f(T)} \le \norm{f}_\omega := \Sum{n \ge 0} \abs{\alpha_n} \omega_n
\end{equation}
and by the submultiplicativity $\omega_{n+m} \le \omega_n \omega_m$ one sees
that $\norm{\cdot}_\omega$ is a function algebra (semi)norm on $\C[z]$.
 
The ``functional calculus'' given by \eqref{intro.e.fc-triv} 
is tailored to the operator $T$
and uses no other information than the growth of the powers of $T$. 
The central question now is: under which conditions 
can one obtain better estimates for $\norm{f(T)}$, i.e., in terms of
weaker function norms?  The conditions we have in mind may involve $T$ (or better:
the semigroup $(T^n)_{n \ge 0}$) or the underlying Banach space. To
recall a famous example: {\em von Neumann's inequality} \cite{vonN51}
states that
if $X= H$ is a Hilbert space and $\norm{T}\le 1$ (i.e., $T$ is a contraction), then 
\begin{equation}\label{intro.e.vonN}
 \norm{f(T)} \le \norm{f}_\infty \qquad \text{for every}\quad f\in \C[z],
\end{equation}
where $\norm{f}_\infty$
is the norm of $f$ in the Banach algebra $\calA = \Ha^\infty(\D)$ of
bounded analytic functions on the open unit disc $\D$.

Von Neumann's result is optimal in the trivial sense that
the estimate \eqref{intro.e.vonN} of course implies that $T$ is a contraction, but also
in the sense that one cannot improve the estimate without further conditions:
If $H= \Ell{2}(\D)$ and $(Th)(z) = zh(z)$ is multiplication
with the complex coordinate, then $\norm{f(T)} = \norm{f}_\infty$ for any
$f\in \C[z]$. A natural question then is to ask 
which operators satisfy the slightly weaker estimate
\[ \norm{f(T)} \lesssim  \norm{f}_\infty \qquad (f\in \C[z])
\]
(called ``polynomial boundedness of $T$''). On a general Banach space
this may fail even for a contraction:
simply take
$X= \ell^1(\Z)$ and $T$ the shift operator, given by 
$(Tx)_n = x_{n+1}$, $n \in \Z$, $x\in \ell^1(\Z)$.
 On the other hand, Lebow \cite{Leb68} has
shown that even on a Hilbert space polynomial boundedness of an operator $T$
may fail if it is only assumed to be {\em power-bounded}, i.e., 
if one has merely $\sup_{n \in \N} \norm{T^n} < \infty$
instead of $\norm{T}\le 1$.  The class of power-bounded operators on Hilbert spaces
is notoriously enigmatic, and it can be considered 
one of the most important problems in operator theory to find
good functional calculus estimates for this class. 

\medskip
Let us shortly comment on the continuous case. Here one is given 
a strongly continuous (in short: $C_0$-)semigroup $(T(s))_{s\ge 0}$ of operators
on a Banach space $X$, and one considers
integrals of the form
\begin{equation}\label{intro.e.sgp-cont}
 \int_{\R_+} T(s) \, \mu(\ud{s}),
\end{equation}
where we assume for simplicity that the support of the measure $\mu$ is bounded.
%Recall that $(T(s))_{s\ge 0}$ is a $C_0$-semigroup if
%$T(0)= \Id$, $T(s+t) = T(s)T(t)$ for $s,t\ge 0$, and for each $x\in X$ the orbit
%\[ \R_+ \pfeil X,\qquad s \tpfeil T(s)x
%\]
%is continuous. 
We shall use only basic results from semigroup theory, and
refer to \cite{ABHN,EN} for further information. The
{\em generator} of the semigroup $(T(s))_{s\ge 0}$ is an, in general unbounded,
closed and desely defined operator $-A$
satisfying
\begin{equation}\label{intro.e.sgp-res}
 (\lambda + A)^{-1} = \int_0^\infty e^{-\lambda s} T(s) \, \ud{s}
\end{equation}
for $\re \lambda$ large enough. The generator is densely defined, i.e.,
its domain $\dom(A)$ is dense in $X$.  In this paper we exclusively
deal with semigroups satisfying a polynomial growth 
$\norm{T(s)} \lesssim (1 + s)^\alpha$ for some $\alpha \ge 0$, and hence
\eqref{intro.e.sgp-res} holds at least for all $\re \lambda > 0$. One writes
$T(s) = e^{-sA}$ for $s\ge 0$ and, more generally, 
\[ (\Lap\mu)(A) := \int_{\R_+} T(s) \,\mu(\ud{s}) 
\]
where 
\[ (\Lap\mu)(z) := \int_{\R_+} e^{-zs} \, \mu(\ud{s})
\]
is the Laplace transform of $\mu$. So in the continuous case
the Laplace transform takes the role of the Taylor series in the discrete case. 
Asking for good estimates for
operators of the form \eqref{intro.e.sgp-cont} is as asking for
{\em functional calculus estimates} for the operator $A$. The 
continuous version of von Neumann's inequality states that
if $X= H$ is a Hilbert space and if $\norm{T(s)} \le 1$ for all $s\ge 0$
(i.e., if $T$ is a {\em contraction semigroup}), then 
\[ \norm{f(A)} \le \norm{f}_\infty \qquad (f = \Lap\mu)
\]
where $\norm{f}_\infty$
is the norm of $f$ in the Banach algebra $\Ha^\infty(\C_+)$ of
bounded analytic functions on the open half place $\C_+ := \{ z \in \C \suchthat
\re z > 0\}$, see \cite[Theorem 7.1.7]{HaaFC}.

There are similarities in the discrete and in the continuous
case, but also characteristic differences. The discrete case is
usually a little more general, shows more irregularities, and
often it is possible to transfer results from the discrete to the continuous
case. (However, this may become quite technical, and we prefer direct
proofs in the continuous case whenever possible.) In the continuous case,
the role of
power-bounded operators is played by
{\em bounded semigroups}, and similar to the discrete case, 
the class of bounded semigroups on Hilbert spaces appears to be rather enigmatic.
In particular, there is a continuous analogue of Lebow's result due to Le Merdy
\cite{LeM00},  cf. also \cite[Section 9.1.3]{HaaFC}.
And there remain some notorious open 
questions involving the functional calculus, e.g., the power-boundedness of
the Cayley transform of the generator, cf. \cite{EisZwa08} and the
references therein.

\medskip

The strongest results in the discrete case obtained so far can be found
in the remarkable paper \cite{Pel82} by Peller from 1982. One of Peller's
results are that  if $T$ is a power-bounded operator on a Hilbert space
$H$, then
\[ \norm{f(T)} \lesssim \norm{f}_{\Be^0_{\infty,1}} \qquad (f\in \C[z])
\]
where is $\Be^0_{\infty,1}(\D)$ is the so-called {\em analytic
Besov algebra} on the disc. (See Section \ref{s.pel} below
for a precise definition). 

In 2005, Vitse \cite{Vit05a} made
a major advance in showing that Peller's Besov class estimate 
still holds true on general Banach spaces if the power-bounded operator $T$
is actually of Tadmor-Ritt type, i.e., satisfies the ``analyticity condition''
\[ \sup_{n \ge 0} \norm{n(T^{n+1} - T^n)} < \infty.
\]
She moreover established in \cite{Vit05b} an analogue for strongly continuous
bounded {\em analytic} semigroups. Whereas Peller's results rest
on Grothendieck's inequality (and hence are particular to Hilbert spaces)
Vitse's approach is based on repeated summation/integration by parts,
possible because of the analyticity assumption.

\medskip

In the present paper we shall complement Vitse's result by devising
an entirely new approach, using our
transference methods (Sections \ref{s.ops} and \ref{s.pel}).
In doing so, we avoid Grothendieck's inequality
and reduce the problem to certain Fourier multipliers on vector-valued
function spaces. By Plancherel's identity,
on Hilbert spaces these are convenient to estimate, 
but one can still obtain positive results on $\Ell{p}$-spaces or
on \UMD\ spaces. Our approach works simultaneously in the discrete
and in the continuous case, and hence we do not only recover
Peller's original result (Theorem \ref{pel.t.pel-disc}) but also
establish a perfect continuous analogue 
(Theorem \ref{pel.t.pel-cont}), conjectured in \cite{Vit05b}.
Moreover, we establish an analogue of the Besov-type estimates
for $\Ell{p}$-spaces and for \UMD\ spaces (Theorem \ref{pel.t.pel-cont-Lp}). 
These 
results, however, are less satisfactory since the algebras
of Fourier multipliers on the spaces $\Ell{2}(\R;X)$ and $\Ell{2}(\Z;X)$  
are not  thoroughly understood if $X$ is not a Hilbert space.

\medskip
In Section \ref{s.gamma} 
we show how our transference methods can also be used to
obtain ``$\gamma$-versions'' of the Hilbert space results.
The central notion here is the so-called $\gamma$-boundedness of
an operator family, a strengthening of operator norm boundedness. It is
related to the notion of $R$-boundedness and plays a major role
in Kalton and Weis' work \cite{KalWei04} on 
the $\Ha^\infty$-calculus. The 'philosophy' behind this
theory is that to each  Hilbert space result based on Plancherel's theorem
there is a corresponding Banach space version, when 
operator norm boundedness is replaced by $\gamma$-boundedness.

We give evidence to this philosophy by showing how our transference results 
enables one to prove $\gamma$-versions of functional calculus
estimates on Hilbert spaces. As examples, we recover the $\gamma$-version of a result of 
Boyadzhiev and deLaubenfels,
first proved by Kalton and Weis in \cite{KalWei04} (Theorem \ref{gamma.t.kw}).  
Then we derive $\gamma$-versions of the Besov calculus theorems in both
the discrete and the continuous forms. The simple idea consists of
going back to the original factorization in the transference method,
but exchanging the function spaces on which the Fourier multiplier operators
act from an $\Ell{2}$-space into a $\gamma$-space. This idea
is implicit in the original proof from \cite{KalWei04} 
and has also been employed in a similar fashion
recently by Le Merdy \cite{LeM10}.

\medskip

Finally, in Section \ref{s.si} we discuss consequences of our 
estimates for full functional calculi and singular integrals for discrete
and continuous semigroups.
For instance,  we prove in Theorem \ref{si.t.main}
that if $(T(s))_{s\ge 0}$ is any strongly continuous semigroup on
a \UMD\ space $X$, then for any $0 < a < b$ the
principal value integral
\[ \lim_{\epsilon \searrow 0} \int_{ \epsilon < \abs{s-b} < a} \frac{T(s)x}{s-b}
\]
exists for all  $x\in X$.  
For $C_0$-{\em groups} this is well-known, cf.~\cite{Haa07b}, but
for semigroups which are not groups, this is entirely new.

\bigskip

\noindent
{\em Terminology.}\\
We use the common symbols $\N,\Z,\R,\C$ for the sets of
natural, integer, real and complex numbers.
In our understanding
$0$ is not a natural number, and we write
\[ \Z_+ := \{ n \in \Z \suchthat n \ge 0\} = \N \cup\{0\}\quad
\text{and}\quad
 \R_+ := \{ t \in \R \st t \ge 0\}.
\]
Moreover, $\D := \{ z\in \C \suchthat
\abs{z} < 1\}$ is the open unit disc, 
$\T := \{ z\in \C \suchthat \abs{z} =1\}$ is the torus, and
$\C_+ := \{ \ z\in \C \suchthat \re z > 0\}$ is the open
right half plane.

We use $X,Y,Z$ to denote (complex) Banach spaces, and $A,B,C$ to denote
closed possibly unbounded operators on them. By $\Lin(X)$ we denote
the Banach algebra of all bounded linear operators on the Banach space $X$,
endowed with the ordinary operator norm.
The domain, kernel and range of an operator $A$ are denoted by 
$\dom(A), \ker(A)$ and $\ran(A)$, respectively. 

The Bochner space of equivalence classes of 
$p$-integrable $X$-valued functions is denoted by $\Ell{p}(\R;X)$. 
If $\Omega$ is a locally compact
space, then $\eM(\Omega)$ denotes the space of all bounded regular
Borel measues on $\Omega$. If $\mu \in \eM(\Omega)$ then 
$\supp\mu$ denotes its topological support. If $\Omega \subset \C$
is an open subset of the complex plane, $\Ha^\infty(\Omega)$ denotes
the Banach algebra of bounded holomorphic functions on $\Omega$,
endowed with the supremum norm $\norm{f}_{\Ha^\infty(\Omega)}
= \sup\{ \abs{f(z)} \suchthat z\in \Omega\}$.

We shall need notation and results from Fourier analysis as
collected in \cite[Appendix E]{HaaFC}. In particular, we use the
symbol $\calF$ for the Fourier transform acting on the space
of (possibly vector-valued) tempered distributions on $\R$, where we 
agree that
\[ \calF\mu(t) := \int_\R e^{-ist}\, \mu(\ud{s})
\]
is the Fourier transform of a bounded measure $\mu \in \eM(\R)$. 
A function $m \in \Ell{\infty}(\R)$ is called a {\em bounded} Fourier
multiplier on $\Ell{p}(\R;X)$ if there is a constant $c\ge 0$ such that
\begin{equation}\label{intro.e.fm}
 \norm{m \cdot\Fourier{f}}_p \le c \norm{\Fourier{f}}_p
\end{equation}
holds true for all $f\in \Ell{p}(\R;X) \cap \Ell{1}(\R;X)$. The
smallest $c$ that can be chosen in \eqref{intro.e.fm} is
denoted by $\norm{\cdot}_{\Mlt_{p,X}}$. This turns the 
space $\Mlt_{p,X}(\R)$ of all bounded Fourier multipliers on 
$\Ell{p}(\R;X)$ into a unital Banach algebra. 

A Banach space $X$ is a \UMD\ space, if and only if the 
function $t \mapsto \sgn t$ is a bounded Fourier multiplier on 
$\Ell{2}(\R;X)$. Such spaces are the right ones  to study
singular integrals for vector-valued functions. In particular,
by results of Bourgain, McConnel and Zimmermann, a 
vector-valued version of the classical Mikhlin theorem holds,
see \cite[Appendix E.6]{HaaFC}  as well as Burkholder's article \cite{Bur01}
and the literature cited there. Each Hilbert space is \UMD, and
if $X$ is \UMD, then $\Ell{p}(\Omega,\Sigma,\mu;X)$ is also \UMD\
whenever $1 < p < \infty$ and $(\Omega,\Sigma,\mu)$ is a measure
space.

The Fourier transform of $\mu \in \ell^1(\Z)$ is
\[ \fourier{\mu}(z) = \sum_{n\in \Z} \mu(n) z^n \qquad (z\in \T).
\]
Analogously to the continuous case, we form the algebra
$\Mlt_{p,X}(\T)$ of functions $m \in \Ell{\infty}(\T)$ which 
induce bounded Fourier multiplier operators on $\Ell{p}(\Z;X)$,
endowed with its natural norm. 

Finally, given a set $A$ and two  real-valued functions $f, g: A \pfeil \R$ we write
\[ f(a) \lesssim g(a) \qquad (a\in A)
\]
to abbreviate the statement that there is $c \ge 0$ such that
$f(a) \le c g(a)$ for all $a\in A$.

\vanish{
All these notions and results of Fourier analysis have their analogues
with  $\R$ replaced by $\Z$. In particular, one can introduce
the Fourier transform

Besides the Fourier transform we shall have need to use the
Laplace transform also.  If $\mu$ is a finite measure on $\R_+)$ then 
\[  (\Lap\mu)(z) := \int_{\R_+} e^{-zs} \, \mu(\ud{s}) \quad (\re z \ge 0)
\]
is its {\em Laplace transform}. 
}

\section{Transference Identities}\label{s.ti}

We introduce the basic idea of transference. Let
$G$ be a locally compact group with left Haar measure $\ud{s}$. 
Let $S \subset G$ be a closed sub-{\em semi}group of $G$
and let 
\[ \rmT : S \pfeil \Lin(X)
\]
be a strongly continuous representation of $S$ on a Banach space $X$.
Let $\mu$ be a (scalar) Borel measure on $S$ such that 
\[ \int_S \norm{T(s)}\, \abs{\mu}(\ud{s}) < \infty, 
\]
and let the operator $\rmT_\mu \in \Lin(X)$ be defined by
\begin{equation}\label{ti.e.sgp-ave} 
\rmT_\mu x := \int_S T(s)x\, \mu(\ud{s}) 
\quad\quad (x\in X).
\end{equation}
The aim of transference is an estimate of 
$\norm{\rmT_\mu}$ in terms of a convolution operator involving $\mu$. 
The idea to obtain such an estimate is, in a first step, purely formal.
To illustrate it we shall need some preparation.

For a (measurable) function $\vphi: S \pfeil \C$ we  denote by 
$\vphi \rmT$ the pointwise product
\[ (\vphi \rmT): S \pfeil \Lin(X), \qquad
s \mapsto  \vphi(s)T(s)
\]
and by $\vphi \mu$ the measure
\[ (\vphi \mu)(\ud{s}) = \vphi(s) \mu(\ud{s}).
\]
In the following we  do not distinguish between a function/measure defined
on $S$ and its extension to $G$ by $0$ on $G \ohne S$.  

Our first lemma expresses the fact that
a semigroup representation induces representations
of convolution algebras on $S$.

\begin{lemma}\label{ti.l.hom}
Let $G, S, \rmT, X$ as above and let $\vphi, \psi : S \pfeil \C$ be functions.
%such that $\vphi \ast \psi$ is defined. 
Then, formally, 
\[ (\vphi \rmT) \ast (\psi \rmT) \, = \, (\vphi \ast \psi)\, \rmT.
\]
\end{lemma}

\begin{proof}
Fix $t\in G$. If $s\in G$ is such that $s \notin S \cap tS^{-1}$ then 
$\vphi(s) = 0$ (in case $s\notin S$) or $\psi(s^{-1}t) = 0$ (in case
$s \notin tS^{-1}$). On the other hand, if $s \in S \cap tS^{-1}$ then 
$s, s^{-1}t \in S$ which implies that $t\in S$ and $T(s)T(s^{-1}t) = T(t)$. 
Hence, formally
%%We have, formally
%\[ \big((\vphi T) \ast (\psi T)\big)(t)= 
%\int_G \vphi(s) T(s) \psi(t-s) T(t-s) \,\ud{s}.
%\]
%for $t\in G$. If $S \cap t-S = \leer$, then for $s\in G$ either
%$\vphi(s) = 0$ or $\psi(t-s) = 0$. Hence $(\vphi \ast \psi)(t) =0$
%and we obtain
%\[ \big((\vphi T) \ast (\psi T)\big)(t)= 0 = (\vphi \ast \psi)(t) T(t).
%\]
%If $t\in G$ is such that $S \cap t-S \not=\leer$, then $t\in S$ and 
%$T(s)T(t-s) = T(t)$ whenever $s, t-s \in S$.
%Hence 
\begin{align*}
\big( (\vphi \rmT) & \ast (\psi \rmT) \big)(t) = \int_G (\vphi \rmT)(s) (\psi \rmT)(s^{-1}t) 
\,\ud{s}
\\ & =
\int_G \vphi(s) \psi(s^{-1}t) T(s) T(s^{-1}t) \,\ud{s}
\\ & =
\int_{S \cap tS^{-1}} \vphi(s)\psi(s^{-1}t) T(s)T(s^{-1}t)\, \ud{s}
\\ & =
\int_{S \cap tS^{-1}} \vphi(s)\psi(s^{-1}t) \, \ud{s}\, T(t)
\\ & = 
\int_G \vphi(s) \psi(s^{-1}t)  \,\ud{s}\, T(t)
= \big( (\vphi\ast \psi)\, \rmT\big)(t).
\end{align*}
\end{proof}

For a function $F: G \pfeil X$ and a measure $\mu$ on $G$ let us
abbreviate
\[ \sprod{F}{\mu} := \int_G F(s)\, \mu(\ud{s})
\]
defined in whatever weak sense. We shall stretch  this notation 
to apply to all cases where it is reasonable. For example, 
$\mu$ could be a vector measure with values in $X'$ or in $\Lin(X)$.

The reflection $F^\sim$  of $F$ is defined by 
\[ F^\sim: G \pfeil X,\quad F^\sim(s) := F(s^{-1}).
\]
If $H: G \pfeil \Lin(X)$ is an operator-valued function, we write
$H \ast F$ for the convolution
\[ (H\ast F)(t) := \int_G H(s)F(s^{-1}t)\, \ud{s} 
\quad \quad (t \in G)
\]
as long as this is well-defined. Also
\[ (\mu \ast F)(t) := \int_G F(s^{-1}t)\,\mu(\ud{s}) \quad \quad (t\in G)
\]
if this is well-defined. 
(Actually, as we are to argue purely formally,
 at this stage we do not bother too much about whether all things
are well-defined. Instead, 
we shall establish a formula first and then
explore conditions under which it is meaningful.) 
The next lemma is almost a triviality.

\begin{lemma}\label{ti.l.fub}
Let $H: G \pfeil \Lin(X)$, $F: G \pfeil X$ and $\mu$ a measure on $G$. Then
\[ \sprod{H \ast F}{\mu} = \sprod{H}{\mu \ast F^\sim}
\]
formally. 
\end{lemma}

\begin{proof}
Writing out the brackets into integrals, it is just Fubini's theorem:
\begin{align*}
\sprod{H \ast F}{\mu } & =
\int_G  \int_G H(s)F(s^{-1}t)\, \ud{s} \, \mu(\ud{t})
= 
\int_G  \int_G H(s)F(s^{-1}t) \, \mu(\ud{t}) \, \ud{s}
\\ & =
\int_G  H(s) \int_G F^\sim(t^{-1}s) \, \mu(\ud{t}) \, \ud{s}
\\ & =
\int_G H(s) (\mu \ast F^\sim)(s)\, \ud{s}
=
\sprod{H}{\mu \ast F^\sim}.
\end{align*}
\end{proof}

If we combine Lemmas \ref{ti.l.hom} and \ref{ti.l.fub} we obtain the following.

\begin{prop}\label{ti.p.ti}
Let $S$ be a closed sub-semigroup of $G$ and  let
$\rmT: S \pfeil \Lin(X)$ be a strongly continuous representation.
Let $\vphi, \vpsi : S \pfeil \C$
and let $\mu$ be a measure on $S$. 
Then, writing $\eta := \vphi \ast \psi$, 
\[  \rmT_{\eta \mu} = \sprod{\rmT}{ (\vphi \ast \psi)\mu }
= \sprod{\vphi \rmT}{ \mu \ast (\vpsi \rmT)^\sim}
\]
formally.
\end{prop}

This result can be interpreted as a {\em factorization} of the
operator $\rmT_{\eta \mu}$ as
\begin{equation}\label{ti.e.ti}
 \begin{xy}
\xymatrix{ \Phi(G;X)\ar[r]^{L_\mu} & \Psi(G;X)\ar[d]^P\\
X\ar[r]^{\rmT_{\eta \mu}}\ar[u]^\iota & X}
\end{xy}
\end{equation}
i.e., \quad  $\rmT_{\eta \mu} = P \nach L_\mu \nach \iota$,
where
\begin{itemize}
\item $\iota$ maps $x\in X$ to the weighted orbit
\[ (\iota x)(s)  := \psi(s^{-1})T(s^{-1})x \quad \quad (s\in G);
\] 
\item $L_\mu$ is the convolution operator with $\mu$  
\[  L_\mu(F) := \mu \ast F;
\]
\item $P$ maps an $X$-valued function on $G$ back to an element of $X$
by integrating against $\vphi \rmT$:
\[ PF := \sprod{\vphi \rmT}{F} = \int_G \vphi(t)T(t)F(t)\, \ud{t},
\]
\item 
$\Phi(G;X), \Psi(G;X)$ are function spaces such that
$\iota: X \pfeil \Phi(G;X)$ and 
$P: \Psi(G;X) \pfeil X$ are meaningful and bounded. 
\end{itemize}
We call a factorization of the form \eqref{ti.e.ti} 
a {\em transference identity}.
It induces a {\em transference estimate} 
\begin{equation}\label{ti.e.te}
\norm{\rmT_{\eta \mu}}_{\Lin(X)} \le \norm{P} \norm{L_\mu}_{\Lin(\Phi(G;X); \Psi(G;X))}
 \norm{\iota}.
\end{equation}

\section{Transference Principles for Groups}\label{s.grp} 

In the present section we shall explain that the classical
transference principle of Berkson-Gillespie-Muhly \cite{BerGilMuh89b} 
for uniformly bounded groups 
and the recent one for  general $C_0$-groups \cite{Haa09b} are instances of
the explained technique.

\subsection{Unbounded $C_0$-groups}\label{grp.ss.unb}

We take $G = S = \R$ and let $\rmU = (U(s))_{s\in \R}: \R \pfeil \Lin(X)$ be a strongly continuous
representation on the Banach space $X$. 
Then $\rmU$ is exponentially bounded, i.e.,
its {\em exponential type}
\[ \theta(\rmU) := \inf\big\{ \omega \ge 0 \st \exists\, M\ge 0 : \norm{U(s)} \le
M \ue^{\omega \abs{s}} \,\,(s\in \R)\big\}
\]
is finite. We choose $\alpha > \omega > \theta(\rmU)$ and
take a measure $\mu$ on $\R$ such that
\[ \mu_\omega := \cosh(\omega\ \cdot) \mu \in \eM(\R)
\]
is a finite measure. Then $\rmU_\mu = \int_\R U\, \ud{\mu}$ is well-defined.
It turns out \cite{Haa09b} that one can factorize
\[ \eta := \frac{1}{\cosh(\omega \cdot)} 
= \vphi \ast \psi 
\]
where $\psi = 1/\cosh(\alpha \cdot)$ and  
$\cosh(\omega \cdot)\vphi = O(1)$. We obtain $\mu = \eta \mu_\omega$
and, writing $\mu_\omega$ for $\mu$ in 
Proposition \ref{ti.p.ti}, 
\begin{equation}\label{grp.unb.e.ti}
 \rmU_\mu = \rmU_{\eta \mu_\omega} = \sprod{\vphi \rmU}{ \mu_\omega \ast (\vpsi \rmU)^\sim} 
= P \nach L_{\mu_\omega} \nach  \iota.
\end{equation}
If $-iA$ is the generator of $\rmU$ and $f = \Fourier{\mu}$ is the Fourier
transform of $\mu$, one writes 
\[ f(A) := \rmU_\mu = \int_\R U(s)\, \mu(\ud{s}),
\]
which is well-defined because the Fourier transform is injective.
Applying the transference estimate \eqref{ti.e.te}  
with $\Phi(\R;X) = \Psi(\R;X) := \Ell{p}(\R;X)$ as the function spaces as  
in \cite{Haa09b} leads to the  estimate
\[ \norm{f(A)} \lesssim \frac{1}{2} 
\Big( \norm{f(\cdot + i\omega)}_{\Mlt_{p,X}(\R)}
+ \norm{f(\cdot - i\omega)}_{\Mlt_{p,X}(\R)} \Big),
\] 
where  $\Mlt_{p,X}(\R)$ denotes the space
of all (scalar-valued) bounded Fourier multipliers  
on $\Ell{p}(\R;X)$.
In the case that $X$ is a \UMD\ space one can now 
use the Mikhlin type result for Fourier multipliers on $\Ell{p}(\R;X)$
to obtain a generalization of the Hieber--Pr\"uss theorem 
\cite{HiePru98} to unbounded groups, see \cite[Theorem 3.6]{Haa09b}.

If $p=2$ and $X=H$,
this Fourier multiplier norm
coincides with the sup-norm by Plancherel's theorem, and by the maximum principle
one obtains  the $\Ha^\infty$-estimate
\begin{equation}\label{grp.e.boydel}
 \norm{f(A)} \lesssim \norm{f}_{\Ha^\infty(\mathrm{St}(\omega))},
\end{equation}
where 
\[ \mathrm{St}(\omega) : =  \{ z\in \C \st \abs{\im z } < \omega\}
\]
is the vertical strip of height $2\omega$, symmetric about the real axis.
This result is originally due to Boyadzhiev and De Laubenfels
\cite{BoydeL94} and is closely related to McIntosh's theorem
on  $\Ha^\infty$-calculus  for sectorial operators with 
bounded imaginary powers from \cite{McI86}, see \cite[Corollary 3.7]{Haa09b}
and \cite[Chapter 7]{HaaFC}.

\subsection{Bounded groups: the classical case}\label{grp.ss.bdd}

The classical transference principle, in the form put forward by 
Berkson, Gillespie and Muhly in \cite{BerGilMuh89b} reads as follows: 
{\em Let $G$ be a locally compact amenable group, let $\rmU = (U(s))_{s\in G}$
be a uniformly bounded, strongly continuous representation of $G$
on a Banach space $X$, and let $p \in [1, \infty)$. Then
\[ \norm{ \int_G U(s)\, \mu(\ud{s})}
\le M^2 \norm{L_\mu}_{\Lin(\Ell{p}(G;X))}
\]
for every bounded measure $\mu \in \eM(G)$.} 
(Here $M := \sup_{s\in G} \norm{U(s)}$.)

\smallskip

We shall review its proof in the special case of $G = \R$
(but the general case is analogous using F\o lner's condition, see \cite[p.10]{CoiWei}).
First, fix $n, N > 0$ and suppose that $\supp(\mu) \subset [-N, N]$. 
Then
\[   \eta := \vphi \ast \psi := 
\frac{1}{2n}\car_{[-n,n]} \ast \car_{[-N-n , N+n]} = 1 
\quad \quad \text{on $[-N,N]$}. 
\]  
So $\eta \mu = \mu$; applying  the transference estimate \eqref{ti.e.te} 
with the function space $\Phi(\R;X) = \Psi(\R;X) := \Ell{p}(\R;X)$ together with 
H\"older's inequality yields
\begin{align*}
 \norm{\rmT_\mu} & \le M^2 \norm{\vphi}_{p'} \norm{\psi}_p 
\norm{L_\mu}_{\Lin(\Ell{p}(\R;X))}
\\ &= M^2 (2n)^{\frac{1}{p'} - 1} (2N + 2n)^{\frac{1}{p}} 
\norm{L_\mu}_{\Lin(\Ell{p}(\R;X))}
\\ & = M^2 \left(1+ \frac{N}{n}\right)^{1/p} \norm{L_\mu}_{\Lin(\Ell{p}(\R;X))}. 
\end{align*}
Finally, let $n \to \infty$ and approximate a general $\mu \in \eM(\R)$ by
measures of finite support.

\begin{rem} 
This proof shows a feature to which we pointed already in the Introduction, but
which was not present in the
case of unbounded groups treated above. Here, an additional 
{\em optimization argument} appears which is based on some freedom 
in the choice of the auxiliary functions $\vphi$ and $\psi$.  
Indeed, $\vphi$ and $\psi$ can vary as long as $\mu = (\vphi \ast \psi) \mu$,
which amounts to $\vphi \ast \psi = 1$ on $\supp(\mu)$. 
\end{rem}

\begin{rem}
A transference principle 
for bounded {\em cosine functions} instead of groups 
was for the first time established and applied in \cite{Haa09a}. 
\end{rem}

\section{A transference principle for discrete and continuous operator semigroups}%
\label{s.ops}

In this section we shall apply the transference method
from Section \ref{s.ti} to operator semigroups, i.e.,
strongly continuous representations of the semigroup 
$\R_+$ (continuous case) or $\Z_+$ (discrete case).

\subsection{The continuous case}

Let $\rmT = (T(s))_{s\ge 0}$ be a strongly continuous
(i.e. $C_0$-) one-parameter semigroup on a (non-trivial)
Banach space $X$. By standard semigroup
theory \cite{EN}, $\rmT$ is exponentially bounded, i.e., there
exists $M, \omega\ge 0$ such that $\norm{T(s)} \le Me^{-\omega s}$ for all
$s\ge 0$. We consider complex measures $\mu$ on $\R_+ := [0, \infty)$ such that
\[ \int_0^\infty \norm{T(s)}\, \abs{\mu}(\ud{s}) < \infty.
\]
If $\mu$ is Laplace transformable and if $f = \Lap\mu$ is 
its Laplace(--Stieltjes)  transform
\[ \Lap\mu(z) = \int_0^\infty e^{-zs}\, \mu(\ud{s}),
\]
then we use (similar to the group
case) the abbreviation
\[ f(A) := \rmT_\mu = \int_0^\infty T(s)\, \mu(\ud{s})
\]
where $-A$ is the generator of the semigroup $\rmT$. 
The mapping $f \mapsto f(A)$  is well-defined since the Laplace transform is injective,
and is called the
{\em Hille--Phillips functional calculus} for $A$, see 
\cite[Section 3.3]{HaaFC} and \cite[Chapter XV]{HilPhi}.

\begin{thm}\label{ops.t.cont}
Let $p \in (1, \infty)$. 
Then there is a constant $c_p \ge 0$ such that
\begin{equation}\label{ops.e.cont}
 \norm{\rmT_\mu} \le c_p \, (1 + \log(b/a))\,  M(b)^2\,  \norm{L_\mu}_{\Lin(\Ell{p}(\R;X))}
\end{equation}
whenever the following hypotheses are satisfied:
\begin{aufziii}
\item $\rmT = (T(s))_{s\ge 0}$\quad  is a $C_0$-semigroup on the Banach space $X$;
\item $0 < a < b < \infty$;
\item $M(b) := \sup_{0\le s \le b} \norm{T(s)}$;
\item $\mu \in \eM(\R_+)$\quad  such that  \quad $\supp(\mu) \subset [a,b]$.
\end{aufziii}
\end{thm}

\begin{proof}
Take $\vphi \in \Ell{p'}(0,b)$, $\psi \in \Ell{p}(0,b)$ such that
$\vphi\ast \psi =1$ on 
$[a,b]$, and let $\eta := \vphi \ast \vpsi$. Then $\eta \mu = \mu$ and 
Proposition \ref{ti.p.ti} yields
\[ \rmT_\mu = \rmT_{\eta \mu} = \sprod{\vphi \rmT}{ \mu \ast (\psi \rmT)^\sim}.
\]
H\"older's inequality leads to a norm estimate
\[ \norm{\rmT_\mu} \le M(b)^2 \norm{\vphi}_{p'} \norm{\psi}_p 
\norm{L_\mu}_{\Lin(\Ell{p}(\R;X))}.
\]
Hence, to prove the theorem it suffices to show that
\[ c(a,b) := \inf\{ \norm{\vphi}_{p'} \norm{\psi}_p \, :\, 
\vphi\ast \psi =1\,\,\text{on} \,\,[a,b]\} \le c_p\, \log(1 + (b/a))
\]
with $c_p$ independent of $a, b$ and $p$. This is done in 
Lemma \ref{app.l.cont-opt}.% in Appendix \ref{app.opt}.
\end{proof}

\begin{rems}\label{ops.r.cont}
\begin{aufziii}
\item The conclusion of the theorem is also true in the case 
$p = 1$ or $p=\infty$, but in this case  
\[ \norm{L_\mu}_{\Lin(\Ell{p}(\R;X))} = 
\norm{\mu}_{\eM(\R)}
\]
is just the total variation norm of $\mu$.
And clearly  $\norm{\rmT_\mu} \le M(b) \norm{\mu}_\eM$, which is stronger than
\eqref{ops.e.cont}.
\item In functional calculus terms, \eqref{ops.e.cont} takes the form
\[ \norm{f(A)} \le 
c_p \, (1 + \log(b/a))\,  M(b)^2\,  \norm{f}_{\AM_{p,X}(\C_+)}
\]
where $f = \Lap\mu$ and
\begin{align*}
\AM_{p,X}(\C_+) :=  
   \{ f \in  \Ha^\infty(\C_+) \suchthat
f(i \, \cdot\,) \in \Mlt_{p,X}(\R)  \}
\end{align*}  
is the (scalar) {\em analytic $\Ell{p}(\R;X)$-Fourier multiplier algebra}, endowed with the norm
\[ \norm{f}_{\AM_{p,X}(\C_+)} := \norm{f(i\cdot)}_{\Mlt_{p,X}(\R)}.
\]
\end{aufziii}
\end{rems}

Let us now state a corollary for
semigroups with polynomial growth type.

\begin{cor}\label{ops.c.cont-pol}
Let $p \in (1, \infty)$. Then there is a constant $c_p \ge 0$ such that
the following is true.
If $-A$ generates a $C_0$-semigroup
$\rmT = (T(s))_{s\ge 0}$ on a Banach space $X$
such that there is $M,\alpha\ge 0$ with
\[ \norm{T(s)} \le M (1 + s)^\alpha \quad\quad (s\ge 0),
\]
then
\begin{equation}\label{ops.e.cont-pol}
 \norm{f(A)} \le c_p M^2 (1+b)^{2\alpha} (1 + \log(b/a))
\norm{f}_{\AM_{p,X}(\C_+)}
\end{equation}
for $0 < a < b < \infty$, $f = \Lap\mu$  and  $\mu \in \eM[a,b]$.
\end{cor}

The case that $\alpha = 0$, i.e., the case of a {\em bounded} semigroup,
is particularly important, hence we state it separately.

\begin{cor}\label{ops.c.cont-bdd}
Let $p \in (1, \infty)$. Then there is a constant $c_p \ge 0$ such that
the following is true. If $-A$ generates a uniformly bounded
$C_0$-semigroup $\rmT = (T(s))_{s\ge 0}$ on a Banach space $X$
then, with $M := \sup_{s\ge 0}\norm{T(s)}$, 
\begin{equation}\label{ops.e.cont-bdd}
 \norm{f(A)} \le c_p M^2 (1 + \log(b/a))
\norm{f}_{\AM_{p,X}(\C_+)}
\end{equation}
for $0 < a < b < \infty$, $f= \Lap\mu$ and  $\mu \in \eM[a,b]$.
\end{cor}

\begin{rem}
If $X = H$ is a Hilbert space and $p=2$, by Plancherel's theorem
and the maximum principle, equation \eqref{ops.e.cont-bdd} becomes
\begin{equation}\label{ops.e.cont-bdd-Hil}
 \norm{f(A)} \lesssim M^2 (1 + \log(b/a))
\norm{f}_{\Ha^\infty(\C_+)}
\end{equation}
where $f= \Lap\mu$ is the Laplace--Stieltjes transform of $\mu$.
A similar estimate has been established by Vitse \cite[Lemma 1.5]{Vit05b}
on a general Banach space $X$, but with the semigroup being
holomorphic and bounded on a sector. 
\end{rem}

\subsection{The discrete case}

We now turn to the situation of a discrete operator semigroup, i.e., the
powers of a bounded operator. Let $T \in \Lin(X)$ be a bounded
operator and $\rmT = (T^n)_{n\in \Z_+}$ the corresponding  
semigroup representation. 
If $\mu \in \ell^1(\Z_+)$ is such that 
$\sum_{n=0}^\infty \abs{\mu(n)} \norm{T^n} < \infty$
then  \eqref{ti.e.sgp-ave} takes the form
\[ \rmT_\mu = \Sum[\infty]{n=0} \mu(n) T^n.
\]
Denoting  $\fourier{\mu}(z) := \sum^\infty_{n=0} \mu(n) z^n$ 
for $\abs{z}\le 1$
we also write $\fourier{\mu}(T) := \rmT_\mu$.

\begin{thm}\label{ops.t.disc}
Let $p \in (1, \infty)$. 
Then there is a constant $c_p \ge 0$ such that
\begin{equation}\label{ops.e.disc}
 \norm{\fourier{\mu}(T)} \le c_p \, (1 + \log(b/a))\,  M(b)^2\,  \norm{L_\mu}_{\Lin(\Ell{p}(\Z;X))}
\end{equation}
whenever the following hypotheses are satisfied:
\begin{aufziii}
\item $T$ is a bounded operator on a Banach space $X$;
\item $a,b\in \Z$ with $1\le a \le b$; 
\item $M(b) := \sup_{0 \le n \le b} \norm{T^n}$;
\item $\mu \in \ell^1(\Z_+)$ such that  $\supp(\mu) \subset [a,b]$.
\end{aufziii}
\end{thm}

\begin{proof}
This is completely analogous to the continuous situation. 
Take $\vphi \in \Ell{p'}(\Z_+)$, $\psi \in \Ell{p}(\Z_+)$ such that
$\vphi\ast \psi =1$ on 
$[a,b]$, and let $\eta := \vphi \ast \vpsi$. Then $\eta \mu = \mu$ and 
Proposition \ref{ti.e.ti} yields
\[ \fourier{\mu}(T)= \rmT_\mu = \rmT_{\eta \mu} = \sprod{\vphi \rmT}{ \mu \ast (\psi \rmT)^\sim}.
\]
H\"older's inequality leads to a norm estimate
\[ \norm{\rmT_\mu} \le M^2 \norm{\vphi}_{p'} \norm{\psi}_p 
\norm{L_\mu}_{\Lin(\Ell{p}(\Z;X))}.
\]
So, similar to the continuous case, one is interested in estimating
\[ c(a,b) := 
\inf \big\{ \norm{\vphi}_{p'} \norm{\psi}_p \, :\, 
\vphi \in \Ell{p'}(\Z_+),\, \psi \in \Ell{p}(\Z_+),\,\,
\vphi\ast \psi =1\,\,\text{on} \,\,[a,b]\big\}.
\]
Applying  Lemma \ref{app.l.disc-opt} concludes the proof.
\end{proof}

\begin{rems}\label{ops.r.disc}
\begin{aufziii}
\item 
As in the continuous case, the assertion remains true for $p = 1, \infty$, but
is weaker than the obvious estimate 
$\norm{\fourier{\mu}}\le M(b) \norm{\mu}_{\ell^1}$.
\item If we write $f = \fourier{\mu}$, \eqref{ops.e.disc} takes the form
\[ 
\norm{f(T)} \le c_p \, (1 + \log(b/a))\,  M(b)^2\,  
\norm{f}_{\AM_{p,X}(\D)}.
\]
Here 
\begin{align*}
\AM_{p,X}(\D)  :=
  \{ f & \in \Ha^\infty(\D) \suchthat
f\res{\T} \in \Mlt_{p,X}(\T) \}
\end{align*}
is the (scalar) {\em analytic $\Ell{p}(\Z;X)$-Fourier multiplier algebra}, endowed
with the norm
\[ \norm{f}_{\AM_{p,X}(\D)} = \norm{f\res{\T}}_{\Mlt_{p,X}(\T)}.
\] 
\end{aufziii}
\end{rems}

Similar to the continuous case we state a consequence for operators with 
polynomially growing powers.

\begin{cor}\label{ops.c.disc-pol}
Let $p \in (1, \infty)$. Then there is a constant $c_p \ge 0$ such that
the following is true.
If $T$ is a bounded operator on a Banach space $X$
such that there is $M,\alpha\ge 0$ with
\[ \norm{T^n} \le M (1+n)^\alpha \quad\quad (n\ge 0),
\]
then
\begin{equation}\label{ops.e.disc-pol}
 \norm{f(T)} \le c_p M^2 (1+b)^{2\alpha} (1 + \log(b/a))
\norm{f}_{\AM_{p,X}(\D)}
\end{equation}
for $f = \fourier{\mu}$,
$\mu \in \ell^1([a,b]\cap \Z)$, $a,b\in \Z$ with $1\le a \le b$.
\end{cor}

\begin{rem}\label{ops.r.relevant-estimate}
For the applications to Peller's theorem in the next section
the exact asymptotics of $c(a,b)$ is irrelevant, and
one can obtain an effective estimate with much less effort.
In the continuous case, the  identity
$c(a,b) = c(1,b/a)$ (cf. the proof of Lemma \ref{app.l.cont-opt}) 
already shows that $c(a,b)$ only depends on 
$b/a$. For the special choice of functions
\[ \vphi = \car_{[0,1]},\quad \psi = \car_{[0,b]}
 \]
one has $\norm{\vphi}_{p'} = 1$ and  
$\norm{\psi}_{p} = b^{1/p}$. 
Consequently $c(1,b) \le b^{1/p}$ and symmetrizing yields
\[ \norm{f(A)} \le M(b)^2 
\left(b/a\right)^{1/ \max(p,p')} \, 
\norm{f}_{\AM_{p,X}(\C)}.
\]
In the discrete case take $\eta$ as in the proof of Lemma 
\ref{app.l.disc-opt} and factorize
\[ \fourier{\eta} = \fourier{\vphi} \cdot \fourier{\psi}=
\frac{1 - z^a}{1-z} \cdot \frac{z}{a (1-z)}. 
\]
Then $\mnorm{\vphi \car_{[0,b]}}_{p'}^{p'} = a$ and  
$\mnorm{\psi \car_{[0,b]}}_{p}^{p} = b/a^p$, hence
\[   c(a,b) \le \norm{\vphi \car_{[0,b]}}_{p'}
\norm{\psi \car_{[0,b]}}_{p} = a^{1/p'} b^{1/p} a^{-1} = (b/a)^{1/p}.
\]
Symmetrizing yields the estimate
\[ \norm{f(T)} \le M^2 
\left(b/a\right)^{1/\max(p,p')} \, \norm{f}_{\AM_{p,X}(\D)}
\]
similar to the continuous case.
\end{rem}

\section{Peller's theorems}\label{s.pel}

The results can be used to obtain a new proof of some
 classical results of Peller's about Besov class functional
calculi for bounded Hilbert space operators with
polynomially growing powers from \cite{Pel82}. In providing
the necessary notions we essentially follow Peller's original
work, changing the notation slightly (cf.~also \cite{Vit05a}).

For an integer $n \ge 1$ let
\[ \vphi_n(k) := \begin{cases} 0, & k\le 2^{n-1},\\
\frac{1}{2^{n-1}} \cdot (k-2^{n-1}), & 2^{n-1} \le k \le 2^n\\
\frac{1}{2^{n}} \cdot (2^{n+1}-k), & 2^{n} \le k \le 2^{n+1}\\
0, & 2^{n+1}\le k.
\end{cases}
\]
That is, $\vphi_n$ is supported in $[2^{n-1}, 2^{n+1}]$, 
zero at the endpoints, $\vphi_n(2^n)= 1$ and linear on each of the
intervals $[2^{n-1}, 2^n]$ and $[2^n, 2^{n+1}]$. Let
$\vphi_0 := (1, 1, 0, \dots)$, then
\[ \Sum[\infty]{n=0} \vphi_n = \car_{\Z_+},
\]
the sum being locally finite. For $s\ge 0$ the {\em Besov class}
$\Be^s_{\infty,1}(\D)$ is defined as the class of  
analytic functions $f$ on the unit disc $\D$ satisfying
\[ \norm{f}_{\Be^s_{\infty, 1}} := 
\Sum[\infty]{n=0} 2^{ns} \norm{\fourier{\vphi_n}\ast f}_{\Ha^\infty(\D)}
%\Sum[\infty]{n=0} 
< \infty.
\]
That is, if $f = \sum_{k\ge 0} \alpha_k z^k$, $\alpha := (\alpha_k)_{k\ge 0}$, 
then
\[ \norm{f}_{\Be^s_{\infty, 1}} = 
%\Sum[\infty]{n=0}  
\Sum[\infty]{n=0}2^{ns} \norm{ \fourier{ \vphi_n\alpha}}_{\Ha^\infty(\D)} < \infty.
\]
Following Peller \cite[p.347]{Pel82}, one has
\[ f\in \Be^s_{\infty,1}(\D) \iff
\int_0^1  (1-r)^{m-s-1} \norm{f^{(m)}}_{\Ell{\infty}(r\T)}\, \ud{r} 
< \infty,
\]
where $m$ is an arbitrary integer such that $m > s$.
Since we only consider $s\ge 0$, we have  
\[ \Be^s_{\infty, 1}(\D) \subset \Ha^\infty(\D)
\]
and it is known that $\Be^s_{\infty, 1}(\D)$ is a Banach algebra
in which  the set of polynomials is dense. 
The following is essentially
\cite[p.354, bottom line]{Pel82}; we give a new proof.

\begin{thm}[Peller 1982]\label{pel.t.pel-disc}
There exists a constant $c\ge 0$ such that the following holds:
Let $X$ be a Hilbert space, and let $T \in \Lin(X)$ such that
\[ \norm{T^n} \le M (1 + n)^\alpha \quad \quad (n \ge 0)
\]
with $\alpha \ge 0$ and $M \ge 1$. Then
\[ \norm{f(T)} \le c \, 9^\alpha \, M^2
\,  \norm{f}_{\Be^{2\alpha}_{\infty,1}(\D)}
\]
for every polynomial $f$.
\end{thm}

\begin{proof}
Let $f = \fourier{\nu} = \sum_{k \ge 0} \nu_n z^n$, and $\nu$ has
finite support. If $n \ge 1$, then $\vphi_n\nu$ has support in 
$[2^{n-1}, 2^{n+1}]$, so we can apply Corollary \ref{ops.c.disc-pol}
with $p =2$ to obtain
\[ \norm{\fourier{\vphi_n \nu}(T)}
\le c_2 M^2 (1 + 2^{n+1})^{2\alpha} (1 + \log 4) 
\norm{\fourier{\vphi_n \nu}}_{\AM_{2,X}(\D)}.
\] 
Since $X$ is a Hilbert space, Plancherel's theorem (and standard
Hardy space theory) yields that
$\AM_{2,X}(\D) = \Ha^\infty(\D)$ with equal norms.
Moreover, $1 + 2^{n+1} \le 3\cdot 2^n$, and hence we obtain
\[ \norm{\fourier{\vphi_n \nu}(T)}
\le c_2 9^\alpha M^2 \cdot 2^{n (2\alpha)} 
\norm{\fourier{\vphi_n \nu}}_{\Ha^\infty(\D)}.
\]
Summing up, we arrive at
\begin{align*}
\norm{f(T)} & \le 
\Sum{n\ge 0} \norm{ \fourier{\vphi_n \nu}(T)}
\\ & \le  \abs{\nu_0} + \abs{\nu_1} M 2^\alpha  
+ c_2 9^\alpha M^2  \Sum{n\ge 1} 
 2^{n (2\alpha)} \norm{\fourier{\vphi_n \nu}}_\infty
\\ & \le c\, 9^\alpha\, M^2\, \, \norm{f}_{\Be^{2\alpha}_{\infty,1}(\D)}
\end{align*}
for some constant $c \ge  0$.
\end{proof}

\begin{rem}\label{pel.r.nikolski}
N.~Nikolski has observed that Peller's Theorem \ref{pel.t.pel-disc}
is only interesting if $\alpha\le 1/2$. Indeed, define
\[ \rmA_{\alpha}(\D) := \left\{ f = \Sum{k\ge 0} a_k z^k \suchthat
\norm{f}_{\rmA_\alpha} := \Sum{k\ge 0} \abs{a_k} (1 + k)^\alpha < \infty
\right\}.
\]
Then $\rmA_\alpha(\D)$ is a Banach algebra with respect to the
norm $\norm{\cdot}_{\rmA_\alpha}$, and one has the obvious estimate
\[ \norm{f(T)} \le M \norm{f}_{\rmA_\alpha} \qquad (f\in \rmA_\alpha(\D))
\]
if  $\norm{T^k} \le M (1+k)^\alpha$, $k \in \N$. This is the 'trivial'
functional calculus for $T$ we mentioned in the Introduction, see 
\eqref{intro.e.fc-triv}.
For $f\in \Be^{\alpha+ 1/2}_{\infty,1}(\D)$ we have
\begin{align*}
\norm{f}_{\rmA_\alpha} & = \abs{a_0} + \Sum{k\ge 0} \Sum{2^{k}\le n < 2^{k+1}}
\abs{a_n} (1+ n)^\alpha
\\ & \le \abs{a_0} + \Sum{k\ge 0} 2^{(k+1)\alpha} \Sum{2^k \le n < 2^{k+1}}
\abs{a_n}
\\ & \le \abs{a_0} + \Sum{k\ge 0}  2^{(k+1)\alpha} 2^{k/2} 
\left( \Sum{2^k \le n < 2^{k+1}} \abs{a_n}^2 \right)^{1/2}
\\ & \le
\abs{a_0} + \Sum{k\ge 0}  2^ \alpha 2^{(\alpha +1/2)k} 
\left(\norm{\fourier{\vphi_{k-1}}\ast f}_2 + 
\norm{\fourier{\vphi_{k}}\ast f}_2 +
\norm{\fourier{\vphi_{k+1}}\ast f}_2 \right)  
\\ & \lesssim
\Sum[\infty]{k=0} 2^{(\alpha+ 1/2)k} \norm{\fourier{\vphi_k} \ast f}_\infty
= \norm{f}_{\Be^{\alpha+ 1/2}_{\infty,1}}
 \end{align*} 
by the Cauchy--Schwarz inequality, Plancherel's theorem and the fact that
$\Ha^\infty(\D) \subset \Ha^2(\D)$. 
This shows that 
$\Be^{\alpha + 1/2}(\D)  \subseteq \rmA_\alpha(\D)$. Hence,
if $\alpha\ge 1/2$, then $2\alpha \ge \alpha+ 1/2$, and therefore
$\Be^{2\alpha}_{\infty,1}(\D) \subset \Be^{\alpha + 1/2}_{\infty,1}(\D)
\subset \rmA_\alpha(\D)$, and the Besov calculus is weaker than the
trivial $\rmA_\alpha$-calculus. 

\smallskip

On the other hand, for $\alpha>0$, 
the example
\[ f(z) = \Sum{n=0}^\infty 2^{-2\alpha n} z^{2^n} \in \rmA_{\alpha}(\D) \ohne
\Be^{2\alpha}_{\infty,1}(\D)
\]
shows that $\rmA_\alpha(\D)$ is {\em not} included into
$\Be^{2\alpha}_{\infty,1}(\D)$, and so the Besov calculus does not
cover the trivial calculus. (By a straightforward
argument one obtains the embedding $\rmA_{\alpha}(\D) \subseteq 
\Be^\alpha_{\infty,1}(\D)$.) 
\end{rem}

\subsection{An analogue in the continuous case}\label{pel.ss.cont}

Peller's theorem has an analogue for
continuous one-parameter semigroups. The role of the unit disc $\D$ is
taken by the right half-plane $\C_+$, the power-series representation
of a function on $\D$ is replaced by a Laplace transform representation
of a function on $\C_+$. However, a subtlety appears that is not present
in the discrete case, namely the possibility  (or even necessity)
to consider also dyadic decompositions ``at zero''. This leads to
so-called ``homogeneous'' Besov spaces, but due to the special form
of the estimate \eqref{ops.e.cont-pol} 
we have to treat the decomposition at $0$
different from the decomposition at $\infty$.

To be more precise, consider the 
partition of unity
\[ \vphi_n(s) := \begin{cases} 0, & 0 \le s\le 2^{n-1},\\
\frac{1}{2^{n-1}} \cdot (s -2^{n-1}), & 2^{n-1} \le s \le 2^n\\
\frac{1}{2^{n}} \cdot (2^{n+1}-s), & 2^{n} \le s \le 2^{n+1}\\
0, & 2^{n+1}\le s
\end{cases}
\]
for $n \in \Z$.
Then $\sum_{n\in \Z} \vphi_n = \car_{(0, \infty)}$, 
the sum being locally finite in $(0,\infty)$.
For $s\ge 0$, an analytic function $f: \C_+ \pfeil \C$ 
is in the (mixed-order homogeneous)
Besov space  $\Be^{0,s}_{\infty,1}(\C_+)$ if 
$f(\infty) := \lim_{t\to \infty} f(t)$ exists and
\begin{align*}
\norm{f}_{\Be^{0,s}_{\infty,1}} & := 
\abs{f(\infty)} + \Sum{n < 0} \norm{\Lap\vphi_n \ast f}_{\Ha^\infty(\C_+)}
\\ & \qquad \qquad \qquad +
\Sum{n\ge 0} 2^{ns} \norm{\Lap\vphi_n \ast f}_{\Ha^\infty(\C_+)} < \infty.
\end{align*}
Here $\Lap$ denotes (as before) the Laplace transform
\[ \Lap \vphi(z) := \int_0^\infty e^{-sz}\vphi(s)\, \ud{s} \quad \quad
(\re  z > 0).
\]
Since we 
are dealing with $s\ge 0$ only, it is obvious that
$\Be^{0,s}_{\infty,1}(\C_+) \subset \Ha^\infty(\C_+)$.
Clearly, our definition of $\Be^{0,s}_{\infty,1}(\C_+)$ is a little sloppy, 
and to make it rigorous we would need to
employ the theory of Laplace transforms of  distributions. 
However,
we shall not need that here, because we shall use only functions of the form
$f = \Lap\mu$, where $\mu$ is a bounded measure {\em with compact support}
in $[0, \infty]$. In this case
\[ \Lap\vphi_n \ast f = \Lap\vphi_n \ast \Lap\mu = 
\Lap(\vphi_n \mu)
\]
by a simple computation.

\begin{thm}\label{pel.t.pel-cont}
There is an absolute constant $c\ge 0$ such that the following holds:
Let $X$ be a Hilbert space, and let $-A$ be the generator
of a strongly continuous semigroup $\rmT = (T(s))_{s\in \R_+}$ on $X$
such that
\[ \norm{T(s)} \le M (1 + s)^\alpha \quad \quad (n \ge 0)
\]
with $\alpha \ge 0$ and $M \ge 1$. Then
\[ \norm{f(A)} \le c \, 9^\alpha\, M^2\,\, 
\norm{f}_{\Be^{0, 2\alpha}_{\infty,1}(\C_+)}
\]
for every $f = \Lap\mu$, $\mu$ being a bounded measure on $\R_+$
 of compact support.
\end{thm}

\begin{proof}
The proof is analogous to the proof of Theorem 
\ref{pel.t.pel-disc}. One has
\[ \mu = f(\infty)\delta_0 +
\Sum{n < 0} \vphi_n\mu
+ \Sum{n \ge 0} \vphi_n\mu
\]
where the first series converges in $\eM[0,1]$ and the
second  is actually finite. Hence
\begin{align*}
\norm{f(A)} & \le \abs{f(\infty)} +
\Sum{n\in \Z} \norm{[\Lap (\vphi_n \mu)](A)}
\\ & \lesssim \abs{f(\infty)} +
\Sum{n \in \Z} M^2 (1 + 2^{n+1})^{2\alpha}
\norm{L_{\vphi_n \mu}}_{\Lin(\Ell{2}(\R;X))}
\\ & = 
\abs{f(\infty)} + M^2 \Sum{n \in \Z} (1 + 2^{n+1})^{2\alpha}
\norm{\Lap\vphi_n \ast f}_{\Ha^\infty(\C_+)}
\\ & \lesssim 
\abs{f(\infty)} + 
M^2  \Sum{n < 0} 2^{2\alpha} \norm{\Lap\vphi_n \ast f}_{\Ha^\infty(\C_+)}
\\ & \qquad \qquad + M^2 \Sum{n \ge 0} ( 3 \cdot 2^{n})^{2\alpha}
\norm{\Lap\vphi_n \ast f}_{\Ha^\infty(\C_+)}
\\
& \le M^2 9^\alpha \norm{f}_{\Be^{0, 2\alpha}_{\infty,1}},
\end{align*}
by Plancherel's theorem and Corollary \ref{ops.c.cont-pol}.
\end{proof}

\begin{rem}
The space $\Be^{0,0}_{\infty,1}(\C_+)$
has been considered by Vitse in \cite{Vit05b} 
under the name $\Be^0_{\infty,1}(\C_+)$,
and we refer to that paper for more information. 
In particular, Vitse proves that
$f\in\; \Be^{0,0}_{\infty,1}(\C_+)$ if and only if $f\in \Ha^\infty(\C_+)$ and
\[ \int_0^\infty \sup_{s\in \R} \abs{f'(t + is)}\, \ud{t} < \infty.
\] 
\end{rem}

Let us formulate the special case $\alpha = 0$ as a corollary.

\begin{cor}\label{pel.c.cont-bdd}
There is a constant $c\ge 0$ such that the following is true.
Whenever $-A$ generates a strongly continuous semigroup
$(T(s))_{s\ge 0}$ on a Hilbert space
such that  $\norm{T(s)} \le M$ for all $s\ge 0$, then
\[ \norm{f(A)} \le c M^2 \norm{f}_{\Be^{0,0}_{\infty,1}(\C_+)}
\]
for all $f = \Lap\mu$, $\mu \in \eM(\R_+)$.
\end{cor}

\begin{rems}
\begin{aufziii}
\item 
Vitse \cite[Introduction, p.248]{Vit05b} in a short note suggests to prove
Corollary \ref{pel.c.cont-bdd} by a discretization argument using Peller's
Theorem \ref{pel.t.pel-disc}
for $\alpha = 0$. This is quite plausible, but no details are given
in \cite{Vit05b} and it seems that  
further work is required to make this approach rigorous.

\item (cf.~Remark \ref{ops.r.relevant-estimate})
%Let us conclude this section with a remarkable observation. 
To prove Theorems \ref{pel.t.pel-disc} and \ref{pel.t.pel-cont} 
we did not make full use
of the logarithmic factor $\log(1 + b/a)$ but only of the
fact that it is constant in $n$ if $[a,b] = [2^{n-1}, 2^{n+1}]$. 
However, as Vitse notes in \cite[Remark 4.2]{Vit05b},
the logarithmic factor appears a fortiori; 
indeed, if $\supp\mu \subset [a,b]$
then if we write
\[ \mu = \sum_{n \in \Z} \vphi_n \mu
\]
the number $N= \card\{ n\in \Z \suchthat \vphi_n \mu \not= 0\}$ 
of non-zero terms in the sum is proportional to $\log(1 + b/a)$.
Hence, for the purposes of functional calculus estimates 
neither Lemma \ref{app.l.cont-opt} nor \ref{app.l.disc-opt}
is necessary.

\item (cf.~Remark \ref{pel.r.nikolski}) 
Different to the discrete case, the Besov estimates are
not completely uninteresting in the case $\alpha \ge 1/2$, because
$\alpha$ affects only the decomposition at $\infty$. 
\end{aufziii}
\end{rems}

\subsection{Generalizations for UMD spaces}

Our proof of Peller's theorems use essentially that the underlying
space is a Hilbert space.  Indeed, we have applied Plancherel's theorem
in order to  estimate the Fourier multiplier norm of a function by its 
$\Ell{\infty}$-norm. 
%It appears to be an open problem whether Hilbert
%spaces are actually characterized by the property
%that every generator of a bounded $C_0$-semigroup on the space
%admits a bounded $\Be^{0,0}_{\infty,1}$-calculus
%actually $X$ to be Hilbert space. 
Hence we do not expect Peller's theorem to be valid on 
other Banach spaces without modifications. In the next section below
we shall show that replacing ordinary boundedness
of an operator family  by the so-called
$\gamma$-boundedness, Peller's theorems carry
over to arbitrary Banach spaces. Here we suggest a different path, namely
to replace the algebra $\Ha^\infty(\C_+)$ in the construction of the Besov
space $\Be^{0,s}_{\infty, 1}$ by the analytic multiplier algebra
$\AM_{p,X}(\C_+)$, introduced in Remark \ref{ops.r.cont}, 2).
We restrict ourselves to the continuous case, leaving 
the discrete version to the reader.

To simplify notation, let us abbreviate $\calA_p := \AM_{p,X}(\C_+)$. 
For $s\ge 0$ and $f: \C_+ \pfeil \C$ we say
 $f \in \Be^{0,s}_{1}[\calA_p]$ if $f\in \Ha^\infty(\C_+)$,
$f(\infty):= \lim_{t\to\infty} f(t)$ exists and 
\[ \norm{f}_{\Be^{0,s}_{1}[\calA_p]}
:= \abs{f(\infty)} + \sum_{n< 0} \norm{\Lap\vphi_n \ast f}_{\calA_p}
+ \sum_{n\ge 0} 2^{ns} \norm{\Lap\vphi_n \ast f}_{\calA_p}
< \infty.
\]
Then the following analogue of Theorem \ref{pel.t.pel-cont} holds, with
a similar proof.

\begin{thm}\label{pel.t.pel-cont-Lp}
Let $p\in (1, \infty)$.
Then there is a constant $c_p\ge 0$ such that the following holds:
Let $-A$ be the generator
of a strongly continuous semigroup $\rmT = (T(s))_{s\in \R_+}$ on 
a Banach space $X$
such that
\[ \norm{T(s)} \le M (1 + s)^\alpha \quad \quad (n \ge 0)
\]
with $\alpha \ge 0$ and $M \ge 1$. Then
\[ \norm{f(A)} \le c_p \, 9^\alpha\, M^2\,\, 
\norm{f}_{\Be^{0, 2\alpha}_{1}[\calM_p]}
\]
for every $f = \Lap\mu$, $\mu$ a bounded measure on $\R_+$
of compact support.
\end{thm}

For $X=H$ is a Hilbert space and $p=2$ one is back at 
Theorem \ref{pel.t.pel-cont}.
For special cases of $X$ --- typically if $X$ is an $\Ell{1}$- or a 
$\Ce(K)$-space ---  one has 
$\Be^{0, 0}_{1}[\calM_p] = \eM(\R_+)$. But if $X$ is a
\UMD\ space, one has positive results. To formulate them
let
\[ \Ha^\infty_1(\C_+) := \{ f\in \Ha^\infty(\C_+) \st z f'(z) \in \Ha^\infty(\C_+)\}
\]
be the {\em analytic Mikhlin algebra}. This is  a Banach algebra with respect
to the norm 
\[ \norm{f}_{\Ha^\infty_1} := \sup_{z\in \C_+} \abs{f(z)} + \abs{z f'(z)}.
\]
If $X$ is a \UMD\ space then the vector-valued version
of the Mikhlin theorem \cite[Theorem E.6.2]{HaaFC}
implies that one has
a continuous inclusion 
\[ \Ha^\infty_1(\C_+) \subset \AM_{p,X}(\C_+)
\]
where the embedding constant depends on $p$ and (the \UMD\ constant of) $X$.
If one defines $\Be_1^{0,s}[\Ha^\infty_1]$ analogously
to $\Be_1^{0,s}[\calM_p]$ above, then we obtain the following.

\begin{cor}\label{pel.c.pel-cont-Lp-UMD}
If $X$ is a \UMD\ space, then 
Theorem \ref{pel.t.pel-cont-Lp} is still valid 
when $\AM_{p,X}(\C_+)$ is replaced by $\Ha^\infty_1(\C_+)$ and
the constant $c_p$ is allowed to depend on (the \UMD-constant of) $X$.
\end{cor}

Now, fix $\theta \in (\pi/2,\pi)$ and consider the sector
\[ \sector{\theta} := \{ z\in \C\ohne\{0\}  \suchthat 
\abs{\arg z} < \theta\}.
\]
Then $\Ha^\infty(\sector{\theta}) \subset \Ha^\infty_1(\C_+)$, 
as follows from an application of the Cauchy integral formula, see
\cite[Lemma 8.2.6]{HaaFC}.
Hence, if we define $\Be_{\infty,1}^{0,s}(\sector{\theta})$
by replacing the space $\Ha^\infty(\C_+)$ in the definition of
$\Be^{0,s}_{\infty,1}(\C_+)$ by $\Ha^\infty(\sector{\theta})$
we obtain the following \UMD-version of Peller's theorem.

\begin{cor}\label{pel.c.pel-cont-Lp-UMD-sector}
Let $\theta \in (\pi/2, \pi)$, let $X$ be a \UMD\ space, and let
$p \in (1, \infty)$. Then there is a constant $c = c(\theta, X, p)$
such that the following holds. 
Let $-A$ be the generator
of a strongly continuous semigroup $\rmT = (T(s))_{s\in \R_+}$ on 
$X$ such that
\[ \norm{T(s)} \le M (1 + s)^\alpha \quad \quad (n \ge 0)
\]
with $\alpha \ge 0$ and $M \ge 1$. Then
\[ \norm{f(A)} \le c \, 9^\alpha\, M^2\,\, 
\norm{f}_{\Be^{0, 2\alpha}_{\infty,1}(\sector{\theta})}
\]
for every $f = \Lap\mu$, $\mu$ a bounded measure on $\R_+$
of compact support.
\end{cor}
Note that Theorem \ref{pel.t.pel-cont} 
above simply says that if $X$
is a Hilbert space, one can choose $\theta = \pi/2$ in  Corollary 
\ref{pel.c.pel-cont-Lp-UMD-sector}.

\begin{rem}
It is natural to ask whether $\Be^{0,s}_1[\Ha^\infty_1]$
or $\Be^{0, 2\alpha}_{\infty,1}(\sector{\theta})$ are actually 
Banach algebras. This is probably not true, as the underlying Banach algebras
$\Ha^\infty_1(\C_+)$ and $\Ha^\infty(\sector{\theta})$ are not
invariant under shifting along the imaginary axis, and hence
are not $\Ell{1}(\R)$-convolution modules. Consequently, Corollaries
\ref{pel.c.pel-cont-Lp-UMD} and \ref{pel.c.pel-cont-Lp-UMD-sector}
are highly unsatisfactory from a
functional calculus point of view.
\end{rem}

\section{Generalizations involving $\gamma$-boundedness}\label{s.gamma}

At the end of the previous section we discussed one 
possible generalization of Peller's theorems, involving
still an assumption on the Banach space and a modification
of the Besov algebra, but no additional assumption on the semigroup.
Here we follow a different path, strengthening the requirements
on the semigroups under consideration. Vitse has shown in 
\cite{Vit05a,Vit05b} that the Peller-type results remain  true 
without any restriction on the Banach space if the
semigroup is bounded analytic
 (in the continuous case), or the operator is a Tadmor--Ritt operator
(in the discrete case). (These two situations correspond to each other
in a certain sense, see e.g. \cite[Section 9.2.4]{HaaFC}.)

Our approach here is based on the ground-breaking work of Kalton and Weis
of recent years, involving the concept of $\gamma$-boundedness.  This
is a stronger notion of boundedness of a set of operators between
two Banach spaces. The ``philosophy'' of the Kalton-Weis approach
is that every Hilbert space theorem which rests on Plancherel's theorem
(and no other result specific for Hilbert spaces) can be 
transformed into a theorem on general Banach spaces, when
operator norm boundedness (of operator families) is replaced by 
$\gamma$-boundedness.  

The idea is readily sketched. 
In the proof of Theorem \ref{pel.t.pel-cont} 
we used
the transference identity \eqref{ti.e.ti} with the function space
$\Ell{2}(\R;X)$ and factorized the operator $\rmT_\mu$
over the Fourier multiplier $L_\mu$. If $X$ is a Hilbert space,
the $2$-Fourier multiplier norm of $L_\mu$ is just 
$\norm{\Lap\mu}_\infty$ and this led to the Besov class estimate.
We now replace the function space $\Ell{2}(\R;X)$ by 
the space $\gamma(\R;X)$; in order to make sure that 
the transference identity \eqref{ti.e.ti} remains valid,
we need that the embedding $\iota$ and the projection $P$
from \eqref{ti.e.ti} are well defined. And this is where the 
concept of $\gamma$-boundedness comes in. Once 
we have established the transference
{\em identity}, we can pass to the transference {\em estimate}; and since
$\Ell{\infty}(\R)$ is also the Fourier multiplier algebra
of $\gamma(\R;X)$, we recover the infinity norm as in the
$\Ell{2}(\R;H)$-case from above. 

We shall now pass to more rigorous mathematics,
starting with a (very brief) introduction to the theory
of $\gamma$-spaces. For a deeper account we refer to \cite{vNe09}.

\subsection{$\gamma$-summing and $\gamma$-radonifying operators}

Let $H$ be a Hilbert space and  $X$ a Banach space. An operator
$T: H \pfeil X$ is called {\em $\gamma$-summing} if
\[ \norm{T}_{\gamma} := 
\sup_F  \,  \Exp\left(\norm{\Sum{e\in F} \gamma_e \tensor Te}_X^2\right)^{1/2}
< \infty,
\]
where the supremum is taken over all finite orthonormal systems
$F \subset H$ and $(\gamma_e)_{e\in F}$ is an independent collection
of standard Gaussian random variables on some probability space. 
It can be shown that in this definition it suffices to consider
only finite subsets $F$ of some {\em fixed}  orthonormal basis of $H$. 
We let
\[ \gamma_\infty(H;X) := \{ T: H \pfeil X \suchthat 
\text{$T$ is $\gamma$-summing}\}
\]
the space of $\gamma$-summing operators of $H$ into $X$. This is
a Banach space with respect to the norm $\norm{\cdot}_{\gamma}$.
The closure in $\gamma_\infty(H;X)$ of the space of finite rank
operators is denoted by  $\gamma(H;X)$, and its elements $T \in  
\gamma(H;X)$ are called {\em $\gamma$-radonifying}. By a theorem
of Hoffman-J\o rgensen and Kwapie\'n, if $X$ does not contain
$\ce_0$ then $\gamma(H;X) = \gamma_\infty(H;X)$, see \cite[Thm.~4.3]{vNe09}.

\smallskip
From the definition of the $\gamma$-norm the following important
{\em ideal property} of the $\gamma$-spaces is quite straightforward
\cite[Thm.~6.2]{vNe09}.

\begin{lemma}[Ideal Property]\label{gamma.l.ideal}
Let $Y$ be  another  Banach space and $K$ another Hilbert space, let
$L : X \pfeil Y$ and $R: K \pfeil H$  be bounded linear operators, and let
$T \in \gamma_{\infty}(H;X)$. Then 
\[  LTR \in \gamma_{\infty}(K;Y)\quad \text{and}\quad
\quad \norm{LTR}_{\gamma} \le \norm{L}_{\Lin(X;Y)} \norm{T}_\gamma
\norm{R}_{\Lin(K;H)}.
\]
If $T \in \gamma(H;X)$, then $LTR \in \gamma(K;Y)$.
\end{lemma}

\smallskip
If $g \in H$ we abbreviate $\konj{g} := \sprod{\cdot}{g}$, i.e.,
$g \tpfeil \konj{g}$ is the canonical conjugate-linear bijection 
of $H$ onto its dual $\konj{H}$. 
Every finite rank operator $T: H \pfeil X$ has the form
\[ T = \Sum[n]{j=1} \konj{g_n} \tensor x_j,
\]
and one can view $\gamma(H;X)$ as a completion of the algebraic 
tensor product $\konj{H} \tensor X$ with respect to the $\gamma$-norm.
Since
\[ \norm{\konj{g}\tensor x}_\gamma = \norm{g}_H \norm{x}_X
= \norm{\konj{g}}_{\konj{H}} \norm{x}_X
\]
for every $g\in H$, $x\in X$, the $\gamma$-norm is a cross-norm. Hence
every {\em nuclear} operator $T: H \pfeil X$ is $\gamma$-radonifying 
and $\norm{T}_\gamma \le \norm{T}_{\text{nuc}}$. (Recall  that $T$ is a nuclear
operator if $T = \sum_{n\ge 0} \konj{g_n} \tensor x_n$ for some
$g_n \in H, x_n \in X$ with 
$\sum_{n\ge 0} \norm{g_n}_H \norm{x_n}_X < \infty$.) The following 
application turns out to be quite useful.

\begin{lemma}\label{gamma.l.nuc}
Let $H,X$ as before, and let $(\Omega,\Sigma,\mu)$ be a measure space. 
Suppose that $f: \Omega \pfeil H$ and $g : \Omega \pfeil X$ 
are (strongly) $\mu$-measurable and 
\[ \int_\Omega \norm{f(t)}_H \norm{g(t)}_X \, \mu(\ud{t}) < \infty.
\]
Then $\konj{f} \tensor g \in \Ell{1}(\Omega; \gamma(H;X))$, and 
$T := \int_\Omega \konj{f} \tensor g \,\ud{\mu} \in \gamma(H;X)$
satisfies
\[ Th = \int_\Omega \sprod{h}{f(t)} g(t)\, \mu(\ud{t}) 
\qquad (h \in H)
\]
and
\[ \norm{T}_\gamma \le \int_\Omega \norm{f(t)}_H \norm{g(t)}_X \, \mu(\ud{t}).
\]
\end{lemma}

Suppose that $H = \Ell{2}(\Omega, \Sigma,\mu)$ for some measure space
$(\Omega,\Sigma,\mu)$.  Every function $u\in \Ell{2}(\Omega;X)$
defines an operator $T_u: \Ell{2}(\Omega) \pfeil X$  by integration:
\[ T_u : \Ell{2}(\Omega) \pfeil X,\qquad T_u(h) = \int_\Omega h\cdot u\, 
\ud{\mu}. 
\]
(Actually, one can do this under weaker hypotheses on $u$, but we
shall have no occasion to use the more general version.)
In this context we identify the operator $T_u$ with the function $u$
and write 
\[ u \in \gamma_{(\infty)}(\Omega; X) \quad \text{in place
of}\quad  T_u \in \gamma_{(\infty)}(\Ell{2}(\Omega);X)).
\]

Extending an idea of  \cite[Remark 3.1]{KalWei04} 
we can use Lemma \ref{gamma.l.nuc}
to conclude
that certain vector-valued functions define $\gamma$-radonifying
operators. Note that $a = -\infty$ or $b = \infty$ are allowed in the following
lemma.

\begin{cor}\label{gamma.c.nuc}
Let $(a,b) \subset \R$, let
$u\in  \Wee{1,1}_{\loc}((a,b);X)$ and let $\vphi : (a,b) \pfeil \C$. 
Suppose that 
one of the following two conditions is satisfied:
\begin{aufziii}
\item $\norm{\vphi}_{\Ell{2}(a,b)}  \norm{u(a)}_X < \infty$ \quad and
\quad $\displaystyle 
\int_a^b \norm{\vphi}_{\Ell{2}(s,b)} \norm{u'(s)}_X \ud{s} < \infty$;

\item  $\norm{\vphi}_{\Ell{2}(a,b)}  \norm{u(b)}_X < \infty$ \quad and
$\displaystyle  
\int_a^b \norm{\vphi}_{\Ell{2}(a,s)} \norm{u'(s)}_X \ud{s} < \infty$.
\end{aufziii}
Then $\vphi \cdot u\in \gamma((a,b);X)$ with respective
estimates for $\norm{\vphi \cdot u}_\gamma$.
\end{cor}

\begin{proof}
In case 1) we use the  representation 
$u(t) = u(a) + \int_a^t u'(s)\, \ud{s}$, 
leading to 
\begin{align*}
\vphi \cdot u & = \vphi \tensor u(a) + \int_a^b \car_{(s,b)} \vphi \tensor
u'(s)\, \ud{s}. 
%\\& = \vphi \tensor u(b) + \int_a^b \car_{(a,s)} \vphi \tensor u'(s)\, \ud{s}.
\end{align*}
Then we apply Lemma \ref{gamma.l.nuc}. In case 2) we start with 
$u(t) = u(b) - \int_t^b u'(s)\,\ud{s}$ and proceed similarly.
\end{proof}

The space $\gamma(\Ell{2}(\Omega);X)$ can be viewed
 as space of generalized $X$-valued functions on $\Omega$.
Indeed, if  $\Omega = \R$ with the Lebesgue measure, 
$\gamma_\infty(\Ell{2}(\R);X)$
is a Banach space of $X$-valued tempered distributions. For such distributions
their Fourier transform is coherently defined via its adjoint action:
$\Fourier T := T \nach \Fourier$, and the ideal property mentioned above
shows that $\calF$ 
restricts to almost isometric isomorphisms
of $\gamma_\infty(\Ell{2}(\R);X)$ and $\gamma(\Ell{2}(\R);X)$.
Similarly, 
the multiplication with some function $m \in \Ell{\infty}(\R)$ extends
via adjoint action coherently to $\Lin(\Ell{2}(\R);X)$, and 
the ideal property above yields that $\gamma_\infty(\Ell{2}(\R);X)$ and
$\gamma(\Ell{2}(\R);X)$ are invariant. Furthermore, 
\[ \norm{T \mapsto mT}_{\gamma_\infty \to \gamma_\infty}
= \norm{m}_\infty
\]
for every $m \in \Ell{\infty}(\R)$. Combining these two facts we obtain
that for each $m \in \Ell{\infty}(\R)$
the {\em Fourier multiplier} operator with symbol $m$
\[ F_m(T)  := \Fourier^{-1} (m\Fourier T) \quad \quad (T \in \Lin(\Ell{2}(\R);X))
\]
is bounded on $\gamma_\infty(\Ell{2}(\R);X)$  and $\gamma(\Ell{2}(\R);X)$
with norm estimate
\[ \norm{F_m(T)}_\gamma \le \norm{m}_{\Ell{\infty}(\R)} \norm{T}_\gamma.
\]
Similar remarks apply in the discrete case $\Omega = \Z$. 

\smallskip
An important result in the theory of $\gamma$-radonifying operators
is the {\em multiplier theorem}. Here one considers a bounded operator-valued
function $T: \Omega \pfeil \Lin(X;Y)$ and asks
under what conditions the  multiplier operator
\[ \calM_T: \Ell{2}(\Omega;X) \pfeil \Ell{2}(\Omega;Y),\qquad
\calM_Tf = T(\cdot)f(\cdot)
\]
is bounded for the $\gamma$-norms. To formulate the result, one needs
a new notion.

Let $X,Y$ be Banach spaces. A collection $\calT \subset \Lin(X;Y)$
is said to be {\em $\gamma$-bounded} if  there is a constant $c\ge 0$ such that
\begin{equation}\label{gamma.e.gamma-bdd}
  \Exp\left(\norm{\Sum{T\in \calT'} \gamma_T Tx_T}_X^2 \right)^{1/2}
\,\le\, c\,\,\,  \Exp\left(\norm{\Sum{T\in \calT'} \gamma_T x_T}_X^2 \right)^{1/2}
\end{equation}
for all finite subsets $\calT' \subset \calT$, 
$(x_T)_{T\in \calT'}\subset X$. (Again, $(\gamma_T)_{T\in \calT'}$
is an independent collection of standard Gaussian random variables on
some probability space.) If $\calT$ is $\gamma$-bounded,
the smallest constant $c$ such that \eqref{gamma.e.gamma-bdd} holds, 
is denoted by
$\gamma(\calT)$
and is called the {\em $\gamma$-bound} of $\calT$.
We are now ready to state the result, established
by Kalton and Weis in \cite{KalWei04}.

\begin{thm}[Multiplier theorem]\label{gamma.t.mult}
Let $H = \Ell{2}(\Omega)$ for some measure space $(\Omega,\Sigma,\mu)$, and let
$X,Y$ be Banach spaces. Let $T : \Omega \pfeil \Lin(X;Y)$ be a strongly
$\mu$-measurable mapping such that
\[ \calT := \{ T(\omega) \suchthat \omega\in \Omega\} 
\] 
is $\gamma$-bounded. 
Then the multiplication operator
\[ \calM_T: \Ell{2}(\Omega)\tensor X \pfeil \Ell{2}(\Omega;Y)\qquad
f \tensor x \tpfeil f(\cdot)T(\cdot)x
\]
extends uniquely to a bounded operator
\[ \calM_T: \gamma(\Ell{2}(\Omega);X) \pfeil \gamma_\infty(\Ell{2}(\Omega);Y)
\]
with\quad $\norm{\calM_T S}_\gamma \le\,  \gamma(\calT) \, \norm{S}_\gamma$,
\quad $(S\in \gamma(\Ell{2}(\Omega);X))$.
\end{thm}

It is unknown up to now whether such a multiplier $\calM_T$
always must have its range in the smaller class $\gamma(\Ell{2}(\Omega);Y)$.
%the general feeling is that this need not be the case. 

\subsection{Unbounded $C_0$-groups}\label{gamma.ss.kw}

Let us  return to 
the main theme of this paper.
In Section \ref{grp.ss.unb} we have applied the transference identities to
unbounded $C_0$-groups in Banach spaces. In the case of a Hilbert space
this yielded a proof of the Boyadzhiev--de Laubenfels theorem, i.e., that
every generator of a $C_0$-group on a Hilbert space has bounded
$\Ha^\infty$-calculus on vertical strips, if the strip height exceeds
the exponential type of the group. The analogue of this result
for general Banach spaces but under $\gamma$-boundedness conditions
is due to Kalton and Weis \cite[Thm.~6.8]{KalWei04}. We give a new proof
using our transference techniques.

Recall that the {\em exponential type} of a $C_0$-group on a Banach
space $X$ is
\[ \theta(U) := \inf\{ \omega \ge 0 \st \exists\, M\ge 0 : \norm{U(s)} \le
M \ue^{\omega \abs{s}} \,\,(s\in \R)\}.
\]
Let us call the number
\[ \theta_\gamma(U)  :=
\inf\{ \omega \ge 0 \st \{ e^{-\omega \abs{s}} U(s)\suchthat s\in \R\}
\,\, \text{is $\gamma$-bounded}\}
\]
the {\em exponential $\gamma$-type} of the group $U$. If 
$\theta_\gamma(U) < \infty$ we call $U$ {\em exponentially $\gamma$-bounded}.
The following is the $\gamma$-analogue of the Boyadzhiev de-Laubenfels
theorem, see equation \eqref{grp.e.boydel}.

\begin{thm}[Kalton--Weis]\label{gamma.t.kw}
Let $-iA$ be the generator of a $C_0$-group 
$(U(s))_{s\in \R}$ on a Banach space $X$. Suppose that
$U$ is exponentially $\gamma$-bounded. Then $A$ has a bounded
$\Ha^\infty(\mathrm{St}(\omega))$-calculus for every $\omega > \theta_\gamma(U)$.
\end{thm}

\begin{proof}
Choose $\theta_\gamma(U) < \omega < \alpha$.
By usual approximation techniques \cite[Proof of Theorem 3.6]{Haa09b} 
it suffices to show an estimate
\[ \norm{f(A)} \lesssim \norm{f}_{\Ha^\infty(\mathrm{St}(\omega))}
\]
only for $f = \Fourier{\mu}$ with  $\mu$ a measure such that 
$\mu_\omega \in \eM(\R)$. (Recall from Section \ref{grp.ss.unb} that
$\mu_\omega(\ud{t}) = \cosh(\omega t)\, \mu(\ud{t})$, so that
$f= \Fourier{\mu}$ 
has a bounded holomorphic extension to $\mathrm{St}(\omega)$.)
By the
transference identity \eqref{grp.unb.e.ti} the operator $f(A)$
factorizes as
\[ f(A) = P \nach L_{\mu_\omega} \nach \iota.
\]
Here $L_{\mu_\omega}$ is convolution with $\mu_\omega$, 
\[  \iota x (s) = \frac{1}{\cosh{\alpha s}} U(-s)x \qquad (x\in X, s\in \R)
\]
and 
\[ PF = \int_\R \psi(t) U(t)F(t)\, \ud{t}.
\]
In Section \ref{grp.ss.unb} this factorization was 
considered to go via the space $\Ell{2}(\R;X)$, 
i.e.,
\[ \iota: X \pfeil \Ell{2}(\R;X),\qquad P: \Ell{2}(\R;X) \pfeil X.
\]
However, the exponential $\gamma$-boundedness of $U$ will allow
us  to replace 
the space $\Ell{2}(\R;X)$ by $\gamma(\Ell{2}(\R);X)$. Once this is ensured,
the estimate is immediate, since 
convolution with $\mu_\omega$ is the Fourier multiplier with 
symbol $\Fourier{\mu_\omega}$. We know that this is
bounded on $\gamma(\Ell{2}(\R);X)$ with a norm
not exceeding 
$\norm{\Fourier{\mu_\omega}}_{\Ell{\infty}(\R)}$, which
by elementary computations and the maximum principle can be majorized by
$\norm{\Fourier{\mu}}_{\Ha^\infty(\mathrm{St}(\omega))}$, 
cf.~Section \ref{grp.ss.unb}.

To see that indeed $\iota: X \pfeil \gamma(\Ell{2}(\R);X)$,  we write
\[ (\iota x)(s)
= \frac{1}{\cosh{\alpha s}} U(-s)x 
=   \left(e^{-\omega \abs{s}} U(-s)\right)\,  \left(
\frac{e^{\omega\abs{s}}}{\cosh{\alpha s}} x\right)
\]
and use the  Multiplier Theorem \ref{gamma.t.mult}
to conclude that 
$\iota: X \pfeil \gamma_\infty(\R;X)$ boundedly.
To see that $\ran(\iota) \subset \gamma (\R;X)$ we employ
a density argument.
If $x\in \dom(A)$, write $\iota x = \psi \cdot u$ with 
\[ \psi(s) = \cosh(\alpha s)^{-1} \quad \text{and}\quad u(s) = U(-s)x
\qquad (s\in \R).
\]
Then $u \in \Ce^1(\R;X)$,  $u'(s) = i U(-s)Ax$,
$\psi\in \Ell{2}(\R)$, and 
\[  
\int_0^\infty \norm{\psi}_{\Ell{2}(s,\infty)} \norm{u'(s)}_X \, \ud{s}, \quad
\int_{-\infty}^0 \norm{\psi}_{\Ell{2}(-\infty,s)} \norm{u'(s)}_X \, \ud{s}\,\,
 <  \,\infty
\] 
Hence, 
\begin{align*}
\iota x = \psi \cdot u & =
\psi \tensor x + \int_0^\infty \car_{(s,\infty)} \psi \tensor u'(s) \,\ud{s}
\\ & \quad\qquad  
- \int_{-\infty}^0 \car_{(-\infty,s)} \psi \tensor u'(s) \,\ud{s} 
\quad \in\, \gamma(\R;X)
\end{align*}
by Corollary \ref{gamma.c.nuc}. 
(One has to apply 1) to the part of $\psi u$ on $\R_+$ and 
2) to the part on $\R_-$.) Since $\dom(A)$ is
dense in $X$, we conclude that $\ran(\iota) \subset \gamma(\Ell{2}(\R);X)$
as claimed.

\smallskip
Finally, we show that $P: \gamma(\Ell{2}(\R);X) \pfeil X $ is well-defined.
Clearly 
\[ P = \left(\text{integrate against $e^{\theta\abs{t}}\vphi(t)$} \right)
\nach \left(\text{multiply with $e^{-\theta\abs{t}} U(t)$}\right)
\]
where $\theta_\gamma(U) < \theta < \omega$. We know that $\vphi(t) 
= O(e^{-\omega\abs{t}})$, so by the Multiplier Theorem \ref{gamma.t.mult},
everythings works out fine.
Note that in order to be able to apply the multiplier theorem,
we have to start
already in $\gamma(\Ell{2}(\R);X)$. And this is why we had to ensure
that $\iota$ maps there in the first place.
\end{proof}

\begin{rem}
Independently of us, 
Le Merdy \cite{LeM10} has recently obtained a $\gamma$-version of
the classical transference principle for bounded groups. 
The method is similar to ours, by re-reading the transference
principle with the $\gamma$-space in place of a Bochner space. 
\end{rem}

\subsection{Peller's theorem --- $\gamma$-version, discrete case}
\label{gamma.ss.pel-disc}

We now turn to the extension of Peller's theorems  
(see  Section \ref{s.pel})
from Hilbert spaces to general Banach spaces. 
We begin with the discrete case.

\begin{thm}\label{gamma.t.pel-disc}
There is an absolute constant $c\ge 0$ such that the following holds:
Let $X$ be a Banach space, and let $T \in \Lin(X)$ such that the set
\[ \calT := \{ (1 +n)^{-\alpha}T^n \suchthat n \ge 0\}
\]
is $\gamma$-bounded. Then
\[ \norm{f(T)} \le c \, 9^\alpha \,   \gamma(\calT)^2 \,\, 
\norm{f}_{\Be^{2\alpha}_{\infty,1}(\D)}
\]
for every polynomial $f$.
\end{thm}

The theorem is a consequence of the following lemma, the arguments being
completely anologous to the proof of Theorem \ref{pel.t.pel-disc}.

\begin{lemma}\label{gamma.l.pel-disc}
There is a constant $c\ge 0$ such that
\[ \norm{\fourier{\mu}(T)} \le c (1 + \log(b/a)) M(b) 
\norm{\fourier{\mu}}_{\Ha^\infty(\D)}
\]
whenever the following hypotheses are satisfied:
\begin{aufziii}
\item $T$ is a bounded operator on a Banach space $X$;
\item $a,b\in \Z$ with $1\le a \le b$; 
\item $M(b) := \gamma\{ T^n \suchthat 0\le n \le b\}$;
\item $\mu \in \ell^1(\Z_+)$ such that  $\supp(\mu) \subset [a,b]$.
\end{aufziii}
\end{lemma}

\begin{proof}
This is analogous to Theorem \ref{ops.t.disc}. 
Take $\vphi, \psi \in \Ell{2}(\Z_+)$
such that $\psi \ast \vphi =1$ on $[a,b]$ and $\supp\vphi, \supp\psi \subset 
[0, b]$. Then
\[ \fourier{\mu}(T) = \sprod{\vphi \rmT}{ \mu \ast (\psi \rmT)^\sim}
= P \nach L_\mu \nach \iota,
\]
see \eqref{ti.e.ti}.
Note that only functions of finite support are involved here, so
$\ran(\iota)\subset \Ell{2}(\Z) \tensor X$. Hence we can take $\gamma$-norms
and estimate
\[ \norm{\fourier{\mu}(T)}
\le \norm{P}_{\gamma(\Ell{2}(\Z);X) \to X} 
\norm{L_\mu}_{\gamma \to \gamma} \norm{\iota}_{X \to \gamma(\Ell{2}(\Z);X)}.
\]
Note that
\[ \iota x = \Big( \rmT^\sim \car_{[-b,0]} \Big) \, (\psi^\sim \tensor x)
\]
so the multiplier theorem yields
\[ \norm{\iota x}_\gamma \le M(b) \norm{\psi^\sim \tensor x}_\gamma
= M(b) \norm{\psi}_2 \norm{x}.
\]
Similarly, $P$ can be decomposed as
\[ P = (\text{integrate against $\vphi$}) \nach (\text{multiply with
$\car_{[0,b]} \rmT$})
\]
and hence the multiplier theorem yields
\[ \norm{P}_{\gamma \to X} \le \norm{\vphi}_2 M(b).
\]
Finally note that 
\[ \norm{L_\mu}_{\gamma \to \gamma} = \norm{\fourier{\mu}}_{\Ha^\infty(\D)}
\]
since --- similar to the continuous case --- all bounded measurable
functions on $\T$
define bounded Fourier multipliers on $\gamma(\Ell{2}(\Z);X)$.
Putting the pieces together we obtain
\[ \norm{\fourier{\mu}(T)} \le M(b)^2 \norm{\vphi}_2 \norm{\psi}_2 
\norm{\fourier{\mu}}_{\Ha^\infty(\D)}
\]
and an application of Lemma \ref{app.l.disc-opt} concludes the proof.
\end{proof}

\subsection{Peller's theorem --- $\gamma$-version, continuous case}

We turn to the continuous version(s) of Peller's theorem.

\begin{thm}\label{gamma.t.pel-cont}
There is an absolute constant $c\ge 0$ such that the following holds:
Let $-A$ be the generator 
of a strongly continuous semigroup $\rmT = (T(s))_{s\ge 0}$
on a Banach space $X$. Suppose that $\alpha \ge 0$ is such that
the set
\[ \calT := \{ (1 +s)^{-\alpha}T(s) \suchthat s \ge 0\}
\]
is $\gamma$-bounded. Then
\[ \norm{f(A)} \le c \, 9^\alpha \, \gamma(\calT)^2 \,\, 
\norm{f}_{\Be^{0, 2\alpha}_{\infty,1}(\C_+)}
\]
for every $f = \Lap\mu$, $\mu$ a bounded measure on $\R_+$ of compact support.
\end{thm}

Let us formulate a minor generalization in the special case 
of bounded semigroups.

\begin{cor}\label{gamma.c.pel-cont-bdd}
There is an absolute constant $c\ge 0$ such that the following holds:
Let $-A$ be the generator 
of a strongly continuous semigroup $\rmT = (T(s))_{s\ge 0}$
on a Banach space $X$ such that the set
\[ \calT := \{ T(s) \suchthat s \ge 0\}
\]
is $\gamma$-bounded. Then
\[ \norm{f(A)} \le c \, \gamma(\calT)^2 \,\, 
\norm{f}_{\Be^{0,0}_{\infty,1}(\C_+)}
\]
for every $f = \Lap\mu$, $\mu$ a bounded measure on $\R_+$.
\end{cor}

The proofs are analogous to the proofs in the Hilbert space case, based 
on the following lemma.

\begin{lemma}\label{gamma.l.cont}
There is a constant $c \ge 0$ such that
\begin{equation}
 \norm{f(A)} \le c \, (1 + \log(b/a))\,  M(b)^2\,  \norm{f}_{\Ha^\infty(\C_+)}
\end{equation}
whenever the following hypotheses are satisfied:
\begin{aufziii}
\item $\rmT = (T(s))_{s\ge 0}$ is a $C_0$-semigroup on the Banach space $X$;
\item $0 < a < b < \infty$;
\item $M(b) := \gamma\{ T(s) \suchthat 0\le s \le b\} $;
\item $f = \Lap\mu$, where $\mu \in \eM(\R_+)$ such that  $\supp(\mu) \subset [a,b]$.
\end{aufziii}
\end{lemma}

\begin{proof}
We re-examine the proof of Theorem \ref{ops.t.cont}. Choose
$\vphi, \psi \in \Ell{2}(0,b)$ such that $\vphi \ast \psi = 1$ on $[a,b]$.
Then 
\[ f(A) = \rmT_\mu =  P \nach L_\mu \nach \iota,
\]
where for $x\in X$ and $F: \R \pfeil X$
\[ \iota x = \psi^\sim  \rmT^\sim  x, \quad 
PF = \int_0^b \vphi(t)T(t)F(t)\, \ud{t}.
\]
We claim that $\iota: X \pfeil \gamma(\R;X)$ with 
\[ 
\norm{\iota}_{X\to \gamma} \le M(b) \norm{\psi}_{\Ell{2}(0,b)}.
\]
As in the case of groups, the estimate follows from the multiplier
theorem; and the fact that $\ran(\iota) \subset \gamma(\R;X)$
(and not just $\gamma_\infty(\R;X)$) comes from a density argument.
Indeed, if $x\in \dom(A)$ then $\iota x = \psi^\sim \cdot u$ with 
$u(s) = T(-s)x$ for $s\le 0$. Since $u \in \Ce^1[-b, 0]$ and 
$\psi^\sim \in \Ell{2}(-b,0)$, Corollary \ref{gamma.c.nuc} 
and the ideal property
yield that  $\iota x = \psi^\sim \cdot u \in \gamma((-b,0);X) \subset 
\gamma(\R;X)$. Since $\dom(A)$ is dense in $X$, 
$\ran(\iota) \subset \gamma(\R;X)$, as claimed.

Note that $P$ can be factorized as
\[ P = \left(\text{integrate against $\vphi$}\right) \nach 
\left(\text{multiply with $\car_{(0,b)} \rmT$}\right)
\]
and so
$\norm{P}_{\gamma \pfeil X} \le M(b) \norm{\vphi}_{\Ell{2}(0,b)}$
by the multiplier theorem.  
We combine these results to obtain
\[ \norm{f(A)} \le M(b)^2 \norm{\vphi}_{\Ell{2}(0,b)}
\norm{\psi}_{\Ell{2}(0,b)} \cdot \norm{f}_{\Ha^\infty(\C_+)}
\]
and an application of Lemma \ref{app.l.cont-opt} concludes the proof.
\end{proof}

\section{Singular Integrals and Functional Calculus}\label{s.si}

\subsection{Functional Calculus}

The results of Sections \ref{s.pel} and \ref{s.gamma} 
provided estimates
of the form
\[ \norm{f(A)} \lesssim \norm{f}_{\Be^{0,2\alpha}_{\infty,1}(\sector{\theta})}
\]
under various conditions on the Banach space $X$, the semigroup $T$
or the angle $\theta$. However, to derive
these estimates we required 
$f = \Lap \mu$, $\mu$ some bounded measure of compact support. 
It is certainly natural to ask whether one can extend the results to all 
$f\in \Be^{0,2\alpha}_{\infty,1}(\sector{\theta})$, i.e., to a proper
Besov class functional calculus.

The major problem here is not the norm estimate, but the 
{\em definition} of $f(A)$ in the first place. (If $f = \Lap \mu$
for a measure $\mu$ with compact support, this problem does not occur.)
Of course one could pass to a closure with respect to the Besov norm,
but this yields a too small function class 
in general.  And it does not show how
this definition of $f(A)$ relates with all the others in the literature, 
especially, with the functional calculus for sectorial operators
\cite{HaaFC} and  the one for half-plane type operators 
\cite{Haa06fpre}.

Unfortunately, although these questions appear to have
quite satisfying answers,  a diligent treatment of them
would extend this already long paper beyond
a reasonable size, so we postpone it to a future publication.

\subsection{Singular Integrals for Semigroups}

A usual consequence of transference estimates is the convergence
of certain singular integrals. It has been known for a long time
that if $(U(s))_{s\in \R}$ is a $C_0$-group 
on a \UMD\ space $X$ then the principal value integral
\[ \int_{-1}^1 U(s)x \, \frac{\ud{s}}{s} 
\]
exists  for every $x\in X$. This was the decisive ingredient
in the Dore--Venni theorem and 
in Fattorini's theorem, as was discussed in \cite{Haa07b}. For semigroups,
these proofs fail and this is not surprising as one has to profit
from cancellation effects around $0$ 
in order to have a principal value integral
converging. Our results from Sections \ref{s.ops} and \ref{s.pel}
now imply that if one shifts the singularity away from $0$ then 
the associated singular integral for a semigroup will converge,
under suitable assumptions on the Banach space or the semigroup. For
groups we gave a fairly general statement in \cite[Theorem 4.4]{Haa09b}.

\begin{thm}\label{si.t.main}
Let $(T(s))_{s\ge 0}$ be a $C_0$-semigroup on a \UMD\ Banach space $X$, 
let $0 < a < b$, and let $g \in \BV[b-a, b+a]$ be such that $g(\cdot + b)$
is even. Then the principal value integral
\begin{equation}\label{si.e.si}
 \lim_{\epsilon \searrow 0} \int_{\epsilon a \le \abs{s- b} \le a}
g(s)T(s)x \, \frac{\ud{s}}{s- b}
\end{equation} 
converges for every $x\in X$.
\end{thm}

\begin{proof}
If $x\in \dom(A)$ then $T(\cdot)x$ is continuously differentiable
and since $g$ is even about the singularity $b$, a well-known argument
shows that the limit \eqref{si.e.si} exists. Hence, by density, one only has
to show that $\sup_{0 < \epsilon < 1} \norm{f_\epsilon(A)}  < \infty$.
In order to establish this, 
define $h(x) = g(ax + b)$ and 
\[ f_\epsilon(z) =  \int_{\epsilon a \le \abs{s- b} \le a}
g(s)e^{-sz}\, \frac{\ud{s}}{s- b}\qquad (z\in \C).
\]
We use Theorem \ref{ops.t.cont} to estimate
\[ \norm{f_\epsilon(A)}
\lesssim  \left(1 + \log\left(\frac{b+a}{b-a}\right)\right)
\norm{f_\epsilon(i \cdot) }_{\Mlt_{p,X}} .
\]
%where $\norm{f_\epsilon (i \cdot)}_{\Mlt_p}$ denotes the norm of the
%Fourier multiplier operator on $\Ell{p}(\R;X)$ with symbol
%$t \mapsto f_\epsilon(it)$. 
Now, by a change of variables,  
\[ f_\epsilon(it) = e^{-itb} \int_{\epsilon \le \abs{s} \le 1}
e^{-iats} \frac{h(s)}{s} \, \ud{s}
= e^{-ibt} \calF( \PV-\frac{h_\epsilon}{s})(at)
%= -2i e^{-itb} \int_0^1 \frac{\sin(st)}{s} h_\epsilon(s)\, \ud{s}
\]
where $h_\epsilon = h \car_{\{\epsilon \le \abs{s} \le 1\}}$. 
It is a standard fact from Fourier multiplier theory that the exponential
factor in front and the dilation by $a$ in the argument do not change
Fourier multiplier norms. So one is reduced to estimate
the $\Mlt_{p,X}$-norms of the functions
\[  m_\epsilon := \calF( \PV-\frac{h_\epsilon}{s}),\qquad (0 < \epsilon <1).
\]
By the \UMD-version of
Mikhlin's theorem,  $\norm{m_\epsilon}_{\calM_p}$
can be estimated by its Mikhlin norm, and by 
\cite[Lemma 4.3]{Haa09b} this in turn can be estimated by the $\BV$-norm
of $m_\epsilon$. But since $\BV[-1,1]$ is a Banach algebra, and
the characteristic functions $\car_{\{\epsilon \le \abs{s} \le 1\}}$ have
uniformly bounded $\BV$-norms for $\epsilon \in (0,1)$, we are done. 
\end{proof}

\begin{rem}
The result is also true on a general
Banach space if  $\{ T(s) \suchthat 0 \le s \le 1\}$ is 
$\gamma$-bounded. The proof is analogous, but 
in place of Theorem \ref{ops.t.cont} one has
to employ Lemma \ref{gamma.l.cont}.
\end{rem}


\appendix
\renewcommand{\appendixname}{\relax}

\section{Two Lemmata}\label{app.opt}

We provide two lemmata concerning an optimization
problem for convolutions on the halfline or the positive integers.

\begin{lemma}[H.-Hyt\"onen]\label{app.l.cont-opt}
Let $p \in (1, \infty)$. For $0 < a < b$ let
\[ c(a,b) := 
\inf\{ \norm{\vphi}_{p'} \norm{\psi}_p \, :\, 
\vphi \in \Ell{p'}(0,b),\, \psi \in \Ell{p}(0,b),\,\,
\vphi\ast \psi =1\,\,\text{on} \,\,[a,b]\}.
\]
Then there are constants $D_p,C_p > 0$ such that
\[  D_p (1 + \log(b/a))  \le c(a,b)
\le C_p (1 + \log(b/a)) 
\]
for all $0 < a < b$.
\end{lemma}

\begin{proof}
We fix $p \in (1, \infty)$.
Suppose that $\vphi \in \Ell{p'}(\R_+)$ and $\psi \in \Ell{p}(\R_+)$ with 
$\vphi \ast \psi = 1$ on $[a,b]$. Then, by H\"older's inequality,
\[ 1 = \abs{(\vphi \ast \psi)(a)}  \le \norm{\vphi}_{p'} \norm{\psi}_{p},
\]
which implies $c(a,b) \ge 1$. Secondly, 
\begin{align*}
 \log(b/a) & = \abs{ \int_a^b (\vphi \ast \psi)(t) \frac{\ud{t}}{t} }
\le  \int_a^b \int_0^t \abs{\vphi(t-s)} \abs{\psi(s)} \, \ud{s} \, \frac{\ud{t}}{t}
\\ &
\le  
\int_0^\infty \int_s^\infty  \frac{ \abs{\vphi(t-s)}}{t}\, \ud{t}\, \abs{\psi(s)}\, \ud{s}
= \int_0^\infty \int_0^\infty \frac{ \abs{\vphi(t)} \abs{\psi(s)} }{t+s} \, \ud{t}\, \ud{s}
\\ &
\le \frac{\pi}{\sin(\pi/p)} \norm{\vphi}_{p'} \norm{\psi}_p.
\end{align*}
(This is ``Hilbert's absolute inequality'', see 
\cite[Chapter 5.10]{GarIneq}.)
This yields
\[ c(a,b) \ge \frac{\sin(\pi/p)}{\pi} \log\left(\frac{b}{a}\right).
\]
Taking both we arrive at
\[ 1 \vee \frac{\sin(\pi/p)}{\pi} \log\left(\frac{b}{a}\right) \le c(a,b).
\]
Since $\sin(\pi/p) \not= 0$, one can find $D_p > 0$ such that
\[ D_p(1 + \log(b/a)) \le 1 \vee \frac{\sin(\pi/p)}{\pi} 
\log(b/a)
\]
and the lower estimate is established.

To prove the upper estimate we note first that without loss of generality we may
assume that $a=1$. Indeed, passing from $(\vphi, \psi)$ to
$(a^{1/p'} \vphi(a\cdot), a^{1/p}\psi(a\cdot))$ reduces the $(a,b)$-case to the
$(1, b/a)$-case and shows that $c(a,b) = c(1, b/a)$.
The idea is now to choose $\vphi, \psi$ in such a way that
\[ (\vphi \ast \psi)(t)  = \begin{cases}
t, & t\in [0,1], \\
1, & t\ge 1,
\end{cases}
\]
and cut them after $b$. Taking Laplace transforms, this means
\[ \big[(\Lap\vphi) \cdot (\Lap\psi)\big](z)  = \frac{1 - e^{-z}}{z^2}
\]
for $\re z > 0$. Fix $\theta \in (0,1)$ and write
\[ \frac{1 - e^{-z}}{z^2} = \frac{(1 -e^{-z})^{(1-\theta)}}{z} \cdot \frac{(1-e^{-z})^{\theta}}{z}.
\]
Now, by the binomial series,
\[  \frac{(1 -e^{-z})^\theta}{z} = \sum_{k=0}^\infty \alpha^{(\theta)}_k \frac{e^{-kz}}{z}
=  \sum_{k=0}^\infty \alpha^{(\theta)}_k \Lap( \car_{(k,\infty)})(z),
\]
and writing $\car_{(k,\infty)} = \sum_{j=k}^\infty \car_{(j,j+1)}$ we see that
we can take
\[ \psi =  \sum_{k=0}^\infty 
\sum_{j=k}^\infty \alpha^{(\theta)}_k  \car_{(j,j+1)}
= \sum_{j=0}^\infty 
\left( \sum_{k=0}^j \alpha^{(\theta)}_k \right)  \car_{(j,j+1)}
\]
and likewise
\[ \vphi 
= \sum_{j=0}^\infty 
\left( \sum_{k=0}^j \alpha^{(1-\theta)}_k \right)  \car_{(j,j+1)}.
\]
Let $\beta^{(\theta)}_j = \sum_{k=0}^j \alpha^{(\theta)}_k$. By standard asymptotic analysis
\[   \alpha^{(\theta)}_k = O\left(\frac{1}{k^{1+\theta}}\right) \quad\text{and}\quad
\beta^{(\theta)}_j = O\left(\frac{1}{(1+j)^\theta}\right)
\]
%(see \cite{...}). 
It is clear that
\[ c(1,b) \le \norm{\vphi \car_{(0,b)}}_{p'} \norm{\psi \car_{(0,b)}}_p.
\]
Now, 
\[ 
\norm{\psi \car_{(0,b)}}_p^p = \int_0^b \abs{\psi(t)}^p\, \ud{t}
= \sum_{j=0}^\infty 
(\beta_j^{\theta})^p  \int_0^b \car_{(j,j+1)}(t)\, \ud{t}
\lesssim \sum_{j=0}^\infty (1+j)^{-\theta p} \gamma_{j,b}
\]
with 
\[ \gamma_{j,b} = \begin{cases}
1, & j \le b-1,\\
b-j,  & j \le b \le j+1,\\
0 & b \le j.
\end{cases}
\] 
With $\theta := 1/p$ this yields 
\[ \norm{\psi \car_{(0,b)}}_p^p \le 1 + \sum_{j=1}^{\floor{b} -1}   \int_{j}^{j+1} \frac{dx}{x} 
+ \frac{ b - \floor{b}}{1 + \floor{b}}
\le 2 + \log(\floor{b}) \le 2(1 + \log b)
\]    
Analogously, noting that $1-\theta = 1 - (1/p) = 1/p'$,
\[ \norm{\vphi \car_{(0,b)}}_{p'}^{p'} \lesssim 2(1 + \log b)
\]
which combines to 
\[ c(1,b) \lesssim (1 + \log b)
\]
as was to prove.
\end{proof}

Now we state and prove an anologue in the discrete case.

\begin{lemma}\label{app.l.disc-opt}
Let $p \in (1, \infty)$. For $a,b\in \N$ with $a\le b$ let
\[ c(a,b) := 
\inf\{ \norm{\vphi}_{p'} \norm{\psi}_p \, :\, 
\vphi \in \Ell{p'}(\Z_+),\, \psi \in \Ell{p}(\Z_+),\,\,
\vphi\ast \psi =1\,\,\text{on} \,\,[a,b]\}.
\]
Then there are constants $C_p,D_p > 0$ such that
\[  D_p (1 + \log(b/a))  \le c(a,b)
\le C_p (1 + \log(b/a)) 
\]
for all $0 < a < b$.
\end{lemma}

\begin{proof}
The proof is similar to the proof of Lemma \ref{app.l.cont-opt}. 
The lower estimate
is obtained in a totally analogous fashion, making use of the discrete
version of Hilbert's absolute inequality 
\cite[Thm. 5.10.2]{GarIneq} and the estimate
\[ \sum_{n=a}^b \frac{1}{n+1} \ge \frac{1}{2} \log(b/a).
\]
For the  upper estimate we let 
\[ \eta(j) := \begin{cases}
j/a, & j=0,1, \dots, a\\
1,   & j\ge a+1,
\end{cases}
\]
and look for a factorization 
$\vphi \ast \psi = \eta$. Considering the Fourier transform
we find
\[ \fourier{\eta}(z) = \frac{z}{a} \frac{1 - z^a}{(1-z)^2}
\]
and so we try (as in the continuous case) the ``Ansatz''
\[ \psi = \frac{z}{a^{\theta}} \frac{(1 - z^a)^{\theta}}{1-z} \quad \text{and}
\quad
\vphi =  \frac{1}{a^{1-\theta}} \frac{(1 - z^a)^{1-\theta}}{1-z}
\]
for $\theta := 1/p$. Note that
\[ 
\psi(z) = \frac{z}{a^\theta (1-z)} \sum_{j=0}^\infty \alpha^{(\theta)}_j z^{aj}
= \frac{z}{a^\theta (1-z)} \sum_{k=0}^\infty \gamma_k z^k,
\]
where 
\[ \gamma_k = \gamma_k(a,\theta) = \begin{cases}
\alpha^{(\theta)}_{k/a} & \text{if}\, a | k\\
0 & \text{else}.
\end{cases}
\]
Consequently,
\[ \psi(z) = \frac{z}{a^\theta} \sum_{n=0}^\infty \left(\sum_{k=1}^n \gamma_k\right) z^n
= \frac{z}{a^\theta} \sum_{n=0}^\infty \beta^{(\theta)}_{\floor{n/a}} z^n
\]
and, likewise,
\[ \vphi(z) = \frac{1}{a^{1-\theta}} \sum_{n=0}^\infty \beta^{(1-\theta)}_{\floor{n/a}} z^n.
\]
As in the continuous case, it suffices to cut off $\vphi$ and $\psi$ after $b$, so
\[ c(a,b) \le 
\norm{\vphi \car_{[0,b]}}_{p'} \norm{\psi \car_{[0,b]}}_p.
\]
Now write $b = ka + r$ with $0\le r < a$ and $k := \floor{b/a}$; then
\begin{align*}
\norm{\psi\car_{[0,b]}}_p^p & \le 
\frac{1}{a} \sum_{n=0}^b  \abs{\beta^{(\theta)}_{\floor{n/a}}}^p
\lesssim
\frac{1}{a} \sum_{n=0}^b  (1 + \floor{n/a})^{-1}
\\ & = \frac{1}{a} \left( \frac{a}{1} + \frac{a}{2} + \dots + \frac{a}{k} +
\frac{r}{k+1}\right) 
\le \sum_{j=1}^{k+1} \frac{1}{j} 
\\ & \le 1 + \int_1^{k+1} \frac{dx}{x}
= 1 + \log(k+1) \le 2(1 + \log(b/a)).
\end{align*}
A similar estimate holds for $\norm{\vphi\car_{[0,b]}}_p^p$.
\end{proof}


%% References


\bigskip

\noindent
{\bf Acknowledgements}\ The work on this paper has occupied me for more
than 2 years now, in which I had discussions with various colleagues and 
friends, a list too long to be given here.
I am particular indebted to N.~Nikolski (Bordeaux) for 
some very valuable remarks about Peller's theorem 
and the analytic Besov classes.
To T.~Hyt\"onen (Helsinki) I owe the proof
of Lemma \ref{app.l.cont-opt} in the appendix; this was a major motivation to continue,
although it eventually turned out that for functional calculus
estimates one can do without it. Finally, I am grateful to my colleagues
in the Analysis Group of the Delft Institute of Applied Mathematics, for
the excellent atmosphere they create.

%\bibliographystyle{alpha}
%\bibliography{../../library/articles,../../library/books,peller}

\begin{thebibliography}{ABHN01}

\bibitem[ABHN01]{ABHN}
Wolfgang Arendt, Charles~J.K. Batty, Matthias Hieber, and Frank Neubrander.
\newblock {\em {Vector-Valued Laplace Transforms and Cauchy Problems.}}
\newblock {Monographs in Mathematics. 96. Basel: Birkh\"auser. xi, 523 p.},
  2001.

\bibitem[Bd94]{BoydeL94}
Khristo Boyadzhiev and Ralph deLaubenfels.
\newblock Spectral theorem for unbounded strongly continuous groups on a
  {H}ilbert space.
\newblock {\em Proc. Amer. Math. Soc.}, 120(1):127--136, 1994.

\bibitem[BGM89]{BerGilMuh89b}
Earl Berkson, T.~Alistair Gillespie, and Paul~S. Muhly.
\newblock Generalized analyticity in {UMD} spaces.
\newblock {\em Ark. Mat.}, 27(1):1--14, 1989.

\bibitem[Blo00]{Blo00}
Gordon Blower.
\newblock Maximal functions and transference for groups of operators.
\newblock {\em Proc. Edinburgh Math. Soc. (2)}, 43(1):57--71, 2000.

\bibitem[Bur01]{Bur01}
Donald~L. Burkholder.
\newblock Martingales and singular integrals in {B}anach spaces.
\newblock In {\em Handbook of the geometry of {B}anach spaces, {V}ol. {I}},
  pages 233--269. North-Holland, Amsterdam, 2001.

\bibitem[Cal68]{Cal68}
Alaberto~P. Calder{\'o}n.
\newblock Ergodic theory and translation-invariant operators.
\newblock {\em Proc. Nat. Acad. Sci. U.S.A.}, 59:349--353, 1968.

\bibitem[CW76]{CoiWei}
Ronald~R. Coifman and Guido Weiss.
\newblock {\em Transference methods in analysis}.
\newblock American Mathematical Society, Providence, R.I., 1976.
\newblock Conference Board of the Mathematical Sciences Regional Conference
  Series in Mathematics, No. 31.

\bibitem[CW77]{CoiWei77}
Ronald~R. Coifman and Guido Weiss.
\newblock Some examples of transference methods in harmonic analysis.
\newblock In {\em Symposia Mathematica, Vol. XXII (Convegno sull'Analisi
  Armonica e Spazi di Funzioni su Gruppi Localmente Compatti, INDAM, Rome,
  1976)}, pages 33--45. Academic Press, London, 1977.

\bibitem[EN00]{EN}
Klaus-Jochen Engel and Rainer Nagel.
\newblock {\em {One-Parameter Semigroups for Linear Evolution Equations.}}
\newblock {Graduate Texts in Mathematics. 194. Berlin: Springer. xxi, 586 p.},
  2000.

\bibitem[EZ08]{EisZwa08}
Tanja Eisner and Hans Zwart.
\newblock The growth of a {$C\sb 0$}-semigroup characterised by its
  cogenerator.
\newblock {\em J. Evol. Equ.}, 8(4):749--764, 2008.

\bibitem[Gar07]{GarIneq}
D.~J.~H. Garling.
\newblock {\em Inequalities: a journey into linear analysis}.
\newblock Cambridge University Press, Cambridge, 2007.

\bibitem[Haa06a]{HaaFC}
Markus Haase.
\newblock {\em The Functional Calculus for Sectorial Operators}.
\newblock Number 169 in Operator Theory: Advances and Applications.
  Birkh\"auser-Verlag, Basel, 2006.

\bibitem[Haa06b]{Haa06fpre}
Markus Haase.
\newblock Semigroup theory via functional calculus.
\newblock Preprint, 2006.

\bibitem[Haa07]{Haa07b}
Markus Haase.
\newblock Functional calculus for groups and applications to evolution
  equations.
\newblock {\em J. Evol. Equ.}, 11:529--554, 2007.

\bibitem[Haa09a]{Haa09a}
Markus Haase.
\newblock The group reduction for bounded cosine functions on {UMD} spaces.
\newblock {\em Math. Z.}, 262(2):281--299, 2009.

\bibitem[Haa09b]{Haa09b}
Markus Haase.
\newblock A transference principle for general groups and functional calculus
  on {UMD} spaces.
\newblock {\em Math. Ann.}, 345:245--265, 2009.

\bibitem[HP74]{HilPhi}
Einar Hille and Ralph~S. Phillips.
\newblock {\em {Functional Analysis and Semi-Groups. 3rd printing of rev. ed.
  of 1957.}}
\newblock {American Mathematical Society, Colloquium Publications, Vol. XXXI.
  Providence, Rhode Island: The American Mathematical Society. XII, 808 p. },
  1974.

\bibitem[HP98]{HiePru98}
Matthias Hieber and Jan Pr{\"u}ss.
\newblock Functional calculi for linear operators in vector-valued {$L\sp
  p$}-spaces via the transference principle.
\newblock {\em Adv. Differential Equations}, 3(6):847--872, 1998.

\bibitem[KW04]{KalWei04}
Nigel Kalton and Lutz Weis.
\newblock The {${H}^\infty$}-functional caluclus and square function estimates.
\newblock unpublished manuscript, 2004.

\bibitem[Leb68]{Leb68}
Arnold Lebow.
\newblock A power-bounded operator that is not polynomially bounded.
\newblock {\em Michigan Math. J.}, 15:397--399, 1968.

\bibitem[LM00]{LeM00}
Christian Le~Merdy.
\newblock A bounded compact semigroup on {H}ilbert space not similar to a
  contraction one.
\newblock In {\em Semigroups of operators: theory and applications (Newport
  Beach, CA, 1998)}, pages 213--216. Birkh\"auser, Basel, 2000.

\bibitem[LM10]{LeM10}
Christian Le~Merdy.
\newblock $\gamma$-{B}ounded representations of amenable groups.
\newblock {\em Adv. Math.}, 224(4):1641--1671, 2010.

\bibitem[McI86]{McI86}
Alan McIntosh.
\newblock Operators which have an ${H}\sb \infty$ functional calculus.
\newblock In {\em Miniconference on operator theory and partial differential
  equations (North Ryde, 1986)}, pages 210--231. Austral. Nat. Univ., Canberra,
  1986.

\bibitem[vN10]{vNe09}
Jan van~Neerven.
\newblock $\gamma$-{R}adonifying {O}perators --- {A} {S}urvey.
\newblock to appear in: Proceedings of the CMA, 2010.


\bibitem[vN51]{vonN51}
Johann von~Neumann.
\newblock Eine {S}pektraltheorie f\"ur allgemeine {O}peratoren eines unit\"aren
  {R}aumes.
\newblock {\em Math. Nachr.}, 4:258--281, 1951.



\bibitem[Pel82]{Pel82}
Vladimir~V. Peller.
\newblock Estimates of functions of power bounded operators on {H}ilbert
  spaces.
\newblock {\em J. Operator Theory}, 7(2):341--372, 1982.

\bibitem[Vit05a]{Vit05a}
Pascale Vitse.
\newblock A band limited and {B}esov class functional calculus for
  {T}admor-{R}itt operators.
\newblock {\em Arch. Math. (Basel)}, 85(4):374--385, 2005.

\bibitem[Vit05b]{Vit05b}
Pascale Vitse.
\newblock A {B}esov class functional calculus for bounded holomorphic
  semigroups.
\newblock {\em J. Funct. Anal.}, 228(2):245--269, 2005.


\end{thebibliography}

\def\cprime{$'$} \def\cprime{$'$} \def\cprime{$'$} \def\cprime{$'$}

\end{document}